\documentclass[11pt]{amsart}
\usepackage{latexsym, graphics, psfrag, amscd, amssymb, graphicx}
\newcommand{\cA}{\mathcal{A}}

\newcommand{\cC}{\mathcal{C}}

\newcommand{\cG}{\mathcal{G}}

\newcommand{\cM}{\mathcal{M}}

\newcommand{\cN}{\mathcal{N}}

\newcommand{\bC}{\mathbb{C}}

\newcommand{\bQ}{\mathbb{Q}}

\newcommand{\bR}{\mathbb{R}}

\newcommand{\bZ}{\mathbb{Z}}

\newcommand{\RP}{\mathbb{RP}}

\newcommand{\ad}{\mathrm{ad}}

\newcommand{\Ad}{\mathrm{Ad}}

\newcommand{\Hom}{\mathrm{Hom}}

\newcommand{\Tr}{\mathrm{Tr}}

\newcommand{\Aut}{\mathrm{Aut}}

\newcommand{\rank}{\mathrm{rank}}

\newcommand{\Prin}{\mathrm{Prin}}

\newcommand{\diag}{\mathrm{diag}}

\newcommand{\Ker}{\mathrm{Ker}}

\newcommand{\fg}{\mathfrak{g}}

\newcommand{\fm}{\mathfrak{m}}

\newcommand{\fr}{\mathfrak{r}}

\newcommand{\fu}{\mathfrak{u}}

\newcommand{\tSi}{\tilde{\Si}}

\newcommand{\tG}{\tilde{G}}

\newcommand{\tP}{\tilde{P}}

\newcommand{\tV}{\tilde{V}}

\newcommand{\ta}{\tilde{a}}

\newcommand{\tb}{\tilde{b}}

\newcommand{\tc}{\tilde{c}}

\newcommand{\te}{\tilde{e}}

\newcommand{\bc}{\bar{c}}

\newcommand{\bd}{\bar{d}}

\newcommand{\bX}{\bar{X}}

\newcommand{\bV}{\bar{V}}

\newcommand{\bv}{\bar{v}}

\newcommand{\Si}{\Sigma}

\newcommand{\ep}{\epsilon}

\newcommand{\ab}{a_1,b_1,\ldots,a_\ell,b_\ell}

\newcommand{\tab}{\ta_1,\tb_1,\ldots,\ta_\ell,\tb_\ell}

\newcommand{\pab}{\prod_{i=1}^\ell[a_i,b_i]}

\newcommand{\hra}{\hookrightarrow}

\newcommand{\define}{\stackrel{\mathrm{def}}{=}}

\newcommand{\fl}[2]{ X_{\mathrm{flat}}^{ {#1}, {#2}}(G) }

\newcommand{\ym}[2]{ X_{\mathrm{YM} }^{{#1},{#2}}(G) }

\newcommand{\Zfl}[1]{ Z_{\mathrm{flat} }^{\ell, {#1}}(G) }

\newcommand{\Zym}[1]{ Z_{\mathrm{YM} }^{\ell, {#1}}(G) }

\newcommand{\flU}[2]{ X_{\mathrm{flat}}^{ {#1}, {#2}}(U(n)) }

\newcommand{\ymU}[2]{ X_{\mathrm{YM} }^{{#1},{#2}}(U(n)) }

\newcommand{\ZymU}[1]{ Z_{\mathrm{YM} }^{\ell, {#1}}(U(n)) }

\newcommand{\flS}[2]{ X_{\mathrm{flat}}^{ {#1}, {#2}}(U(1)) }

\newtheorem{thm}{Theorem}[section]

\newtheorem{lm}[thm]{Lemma}

\newtheorem{rem}[thm]{Remark}

\newtheorem{cor}[thm]{Corollary}

\newtheorem{pro}[thm]{Proposition}

\newtheorem{ex}[thm]{Example}

\newtheorem{no}[thm]{Notation}

\begin{document}

\parskip=0.35\baselineskip

\baselineskip=1.2\baselineskip

\title{Yang-Mills Connections on Nonorientable Surfaces}
\author{Nan-Kuo Ho}
\address{ Department of
Mathematics, National Cheng-Kung University, Taiwan}
\thanks{The first author was partially supported by Natural Sciences and Engineering
Research Council of Canada Postdoctoral Fellowship.}
\email{nankuo@mail.ncku.edu.tw}

\author{Chiu-Chu Melissa Liu}
\address{Department of Mathematics, Northwestern University \and
Department of Mathematics, Columbia University}
\email{ccliu@math.northwestern.edu}

\dedicatory{Dedicated to the memory of Raoul Bott}


\subjclass{53}
\date{\today}

\begin{abstract}
In \cite{ym}, Atiyah and Bott studied Yang-Mills functional over a Riemann
surface from the point of view of Morse theory. We generalize their study
to all closed, compact, connected, possibly nonorientable surfaces.
\end{abstract}

\maketitle

\tableofcontents

\section{Introduction}

Let $G$ be a compact, connected Lie group.
In \cite{ym}, Atiyah and Bott identified the affine space
$\cA$ of connections on a principal $G$-bundle $P$ over a Riemann surface
with the affine space $\cC$ of holomorphic structures on
$P^\bC=P\times_G G^\bC$, where $G^\bC$ is the complexification of $G$.
The identification $\cA \cong \cC$ is an isomorphism of affine spaces, thus
a diffeomorphism.
It was conjectured in \cite{ym} that under this identification
the Morse stratification of the Yang-Mills functional on $\cA$ exists and coincides
with the stratification of $\cC$ from algebraic geometry \cite{HN, Ra}.
The conjecture was proved by Daskalopoulos in \cite{da} (see also \cite{Rad}).
The top stratum $\cC_{ss}$ of $\cC$ consists of semi-stable holomorphic structures
on $P^\bC$. Atiyah and Bott showed  that the stratification of $\cC$ is
$\cG^\bC$-perfect, where $\cG^\bC=\Aut(P^\bC)$.
It has strong implications on
the topology of the moduli space $\cM(P)$ of $S$-equivalence classes of semi-stable
holomorphic structures on $P^\bC$.
When $\cM(P)$ is smooth, Atiyah and Bott found a complete set of
generators of the cohomology groups $H^*(\cM(P);\bQ)$ and recursive
relations which determine the Poincar\'{e} polynomial $P_t(\cM(P);\bQ)$.
When $\cM(P)$ is singular, their results give generators of the equivariant
cohomology groups $H_{\cG^{\bC}}^*(\cC_{ss};\bQ)$ and formula for the equivariant
Poincar\'{e} series $P^{\cG^\bC}_t(\cC_{ss};\bQ)$.

Under the isomorphism $\cA\cong\cC$, the top stratum $\cC_{ss}$ corresponds
to $\cA_{ss}$ which is the stable manifold of $\cN_{ss}$, the set of  central
Yang-Mills connections, where the
Yang-Mills functional achieves its absolute minimum \cite{ym, da}.
When  the absolute minimum is zero, $\cN_{ss}$ is the set of flat connections
(connections with zero curvature). By \cite[Theorem C]{da},
$$
\cM(P)\cong \cN_{ss}/\cG
$$
where $\cG=\Aut(P)$. So $\cM(P)$ can be identified with the moduli space
of gauge equivalence classes of central Yang-Mills connections on $P$.
When the absolute minimum of the Yang-Mills functional is zero, or equivalently,
the obstruction class $o(P)\in H^2(\Si,\pi_1(G))$ is torsion,  $\cM(P)$ is
the moduli space of gauge equivalence classes of flat connections on $P$.
It is known that flat $G$-connections give rise to representations
$\pi_1(\Si)\to G$, where $\pi_1(\Si)$ is the fundamental group of the
base Riemann surface $\Si$ of $P$.  More precisely,
$$
\bigcup_{\tiny \begin{array}{c}P\in \Prin_G(\Si)\\ o(P)\ \textup{torsion}\end{array}}
\cM(P)=\Hom(\pi_1(\Si),G)/G
$$
where $G$ acts on the representation variety $\Hom(\pi_1(\Si),G)$ by
conjugation. Yang-Mills $G$-connections (critical points of the Yang-Mills functional)
give rise to representations $\Gamma_\bR(\Si)\to G$, where $\Gamma_\bR(\Si)$
is the central extension of $\pi_1(\Si)$ \cite[Section 6]{ym}.

In this paper, we study the Yang-Mills functional on the space of connections
on a principal $G$-bundle $P$ over a closed, compact, connected, {\em nonorientable}
surface $\Si$. The pull back $\tP$ of $P$ to the
orientable double cover $\pi:\tSi\to \Si$ is always topologically trivial (Proposition \ref{thm:top_trivial}), and
$A\mapsto \pi^* A$ gives an inclusion from the space $\cA$ of connections on $P$ into
the space $\tilde{\cA}$ of connections on $\tP$. The Yang-Mills functional on
$\cA$ is the restriction of that on $\tilde{\cA}$. In the nonorientable case,
 the absolute minimum of the Yang-Mills functional is always zero, achieved by
flat connections (see for example \cite{HL3}). We have
$$
\bigcup_{P\in \Prin_G(\Si)}\cM(P)=\Hom(\pi_1(\Si),G)/G
$$
where $\cM(P)$ is the moduli space of gauge equivalence classes of flat
connections on $P$.

Let $\Si$ be a compact, connected, nonorientable surface without boundary.
Then $\Si$ is diffeomorphic to the connected sum of $m>0$ copies of $\RP^2$,
and the Euler characteristic $\chi(\Si) =2-m$.
We derive the following results in this paper:
\begin{enumerate}
\item We establish an exact correspondence between the gauge equivalence classes
of Yang-Mills $G$-connections on $\Si$ and conjugacy classes of representations
$\Gamma_\bR(\Si)\to G$, where $\Gamma_\bR(\Si)$ is the {\em super central
extension} of $\pi_1(\Si)$. (Section \ref{sec:YMrep})
\item  We show that the moduli space of gauge equivalence classes of
flat connections on any fixed principal $G$-bundle $P$ over $\Si$
is nonempty and connected if $\chi(\Si)<0$.
This extends \cite[Theorem 5.2]{HL3} to the case
$\Si=4\RP^2$. (Section \ref{sec:components})

\item  When $G=U(n)$, we give an explicit description of the $\cG$-equivariant Morse stratification of
the Yang-Mills functional, compute the Morse index of each stratum,
and relate lower strata to top strata of spaces of $U(m)$-connections ($m<n$)
on $\Si$ and on its orientable double cover. This reduction also gives us a reduction
of equivariant Poincar\'{e} series. (Section \ref{sec:Unonorientable})
\end{enumerate}
We will describe the reduction (3) for other classical groups in a subsequent work
\cite{HL5}.

In the orientable case, the reduction (3) and the understanding
of the topology of the gauge group are sufficient to determine
the equivariant Poincar\'{e} series of the top stratum recursively
(by induction on dimension of the group $G$). In the nonorientable case,
we need to compute the difference of the equivariant Morse and Poincar\'{e} series, which
vanishes in the orientable case due to equivariant perfectness of the stratification.
We will address this in future works.

Using the Morse theory for the Yang-Mills functional
over a closed (orientable or nonorientable) surface
(studied in \cite{ym} and in this paper, respectively),
D. Ramras proved an Atiyah-Segal theorem for surface
groups in \cite{R}: for any closed surface $\Si \neq S^2, \RP^2$,
$K_{\mathrm{def}}^*(\pi_1(\Si))\cong K^*(\Si)$ for
$*>0$ when $\Si$ is orientable, and for $*\geq 0$ when
$\Si$ is nonorientable, where $K_{\mathrm{def}}^*$ are
Carlsson's deformation $K$-groups.

For the purpose of Morse theory we should consider the Sobolev
space of $L_{k-1}^2$ connections  $\cA^{k-1}$ and the group
of $L_k^2$ gauge transformations $\cG^k$ and $(\cG^\bC)^k$,
where $k\geq 2$. We will not emphasize the regularity issues
through out the paper, but refer the reader to \cite[Section 14]{ym}
and \cite{da} for details.

We now give a clear description of the remaining sections.
In Section \ref{sec:flat}, we review various representation varieties
of flat connections, and show that the
pull-back of any principal $G$-bundle over a nonorientable surface to
its orientable double cover is topologically trivial.
In Section  \ref{sec:functional}, we review definitions of the Yang-Mills functional
and Yang-Mills connections over an orientable surface, and give corresponding definitions
for a nonorientable surface. We describe involutions on the principal $G$-bundles and on the
space of connections induced by the deck transformation on the orientable double cover
of the nonorientable surface.
Section \ref{sec:YMrep} contains our main construction and justification.
We introduce the {\em super central extension} of the fundamental
group of a surface; it is the central extension if and only
if the surface is orientable.
We establish a precise correspondence between Yang-Mills connections and
representations of super central extension. We introduce representation varieties
for Yang-Mills connections, and describe an involution on symmetric representation
varieties induced by the deck transformation on the orientable double cover.
We also introduce extended moduli spaces for nonorientable surfaces.
In Section \ref{sec:general}, we discuss the $\cG$-equivariant
Morse stratification and reduction
for general compact connected Lie groups. As a byproduct, we reproduce and extend the results
on connected components of the moduli space of flat connections over closed
nonorientable surfaces \cite{HL2,HL3}.
We specialize to the case $G=U(n)$ in Section \ref{sec:Uorientable} (orientable case)
and Section \ref{sec:Unonorientable} (nonorientable case). We give explicit description
of the $\cG$-equivariant Morse stratification of the space of connections.
The main reference of Section \ref{sec:Uorientable} of this paper is \cite{ym}.
In \cite[Section 7]{ym}, the reduction is derived at the level of strata, which are infinite
dimensional manifolds.
Knowing the exact correspondence between Yang-Mills connections and representations, we work
mainly at the level of representation varieties which are finite dimensional.

\paragraph{\bf Acknowledgments}
This project was motivated by Raoul Bott's, Kenji Fukaya's and Lisa Jeffrey's
questions on betti numbers of moduli spaces of flat connections
over nonorientable surfaces.
We thank Paul Goerss, Victor Guillemin, Frances Kirwan, Lisa Jeffrey,
Jun Li, Eckhard Meinrenken, Mudumbai S. Narasimhan, Daniel Ramras,
Paul Selick, Reyer Sjamaar, Jonathan Weitsman,
Graeme Wilkin, Christopher Woodward, and Siye Wu for helpful conversations.
We thank Lisa Jeffrey, Eckhard Meinrenken, and Christopher Woodward for
comments and suggestions on our draft.

\section{Flat Connections and Representations} \label{sec:flat}

Let $G$ be a compact, connected Lie group, and let $P$ be a principal $G$-bundle on
a closed, compact, connected surface $\Si$. We say a connection $A\in\cA(P)$ is flat if
its curvature vanishes. Let $\cN_0(P)\subset \cA(P)$ be the set
of flat connections on $P$. Note that when $\Si$ is orientable, $\cN_0(P)$ can be empty.

\subsection{Representation varieties of flat connections}\label{sec:Xfl}
We first introduce some notation similar to that in \cite[Section 2.3]{HL3}.
Let $\Si^\ell_0$ be the closed, compact, connected, orientable surface
with $\ell\geq 0$ handles. Let $\Si^\ell_1$ be the connected sum
of $\Si^\ell_0$ and $\RP^2$, and let $\Si^\ell_2$ be the connected
sum of $\Si^\ell_0$ and a Klein bottle. Any closed, compact, connected surface
is of the form $\Si^\ell_i$, where $\ell$ is a nonnegative integer, $i=0,1,2$.
$\Si^\ell_i$ is orientable if and only if $i=0$. Use $1$ as the identity of $\pi_1(\Si)$ and $e$
as the identity of $G$. We have
\begin{eqnarray*}
\pi_1(\Si^\ell_0)
&=&\langle A_1,B_1,\ldots,A_\ell, B_\ell \mid \prod_{i=1}^\ell [A_i,B_i]=1 \rangle\\
\pi_1(\Si^\ell_1)
&=&\langle A_1,B_1,\ldots,A_\ell, B_\ell, C \mid \prod_{i=1}^\ell [A_i,B_i]=C^2 \rangle\\
\pi_1(\Si^\ell_2)
&=&\langle A_1,B_1,\ldots,A_\ell, B_\ell, D,C \mid \prod_{i=1}^\ell [A_i,B_i]=C D C^{-1} D \rangle
\end{eqnarray*}

It is known that a flat connection gives rise to
a homomorphism $\pi_1(\Si)\to G$. Introduce representation varieties
\begin{eqnarray*}
\fl{\ell}{0}&=&\{(\ab) \in G^{2\ell} \mid \prod_{i=1}^\ell [a_i,b_i]=e \} \\
\fl{\ell}{1}&=&\{(\ab,c) \in G^{2\ell+1} \mid \prod_{i=1}^\ell [a_i,b_i]=c^2 \} \\
\fl{\ell}{2}&=&\{(\ab,d,c) \in G^{2\ell+2} \mid \prod_{i=1}^\ell [a_i,b_i]=c d c^{-1} d \}
\end{eqnarray*}
Then
$$
\bigcup_{P\in \Prin_G(\Si^\ell_i)}\cN_0(P)/\cG_0(P)
\cong \Hom(\pi_1(\Si^\ell_i),G)\cong \fl{\ell}{i}
$$
for $i=0,1,2$, $\ell\geq 0$, where $\cG_0(P)$ is the
based gauge group which consists
of gauge transformations on $G$ that take value the identity $e\in G\cong \Aut(P_{x_0})$ at
a fixed point of $x_0\in \Si^\ell_i$. Let $\cG(P)$ be the gauge group. Then
$$
G=\cG(P)/\cG_0(P).
$$

\begin{ex}\label{flS}
\begin{eqnarray*}
\flS{\ell}{0}&=&\{(\ab) \in U(1)^{2\ell}\} = U(1)^{2\ell} \\
\flS{\ell}{1}&=&\{(\ab,c) \in U(1)^{2\ell+1} \mid c^2=1 \}= U(1)^{2\ell}\times \{\pm 1\} \\
\flS{\ell}{2}&=&\{(\ab,d,c) \in U(1)^{2\ell+2} \mid d^2=1 \}= U(1)^{2\ell+1}\times \{\pm 1\}
\end{eqnarray*}
\end{ex}

 Let $G$ act on $G^{2\ell+i}$ by
$$
g\cdot(c_1,\ldots, c_{2\ell+i})= (g c_1 g^{-1},\ldots, g c_{2\ell+i} g^{-1}).
$$
This action preserves the subset $\fl{\ell}{i}\subset G^{2\ell+i}$, so
$G$ acts on $\fl{\ell}{i}$, and
$$
\bigcup_{P\in \Prin_G(\Si^\ell_i)}\cN_0(P)/\cG(P)
\cong \Hom(\pi_1(\Si^\ell_i),G)/G\cong \fl{\ell}{i}/G
$$
is the moduli space of gauge equivalence classes of flat $G$-connections
on $\Si^\ell_i$.

\begin{no}
In the rest of this paper, we will use the following notation:

Denote the $2\ell$-vector by $V=(\ab)\in  G^{2\ell}$. Define $\fm(V)$ and $\fr(V)$ by
\begin{equation} \label{eqn:commutator}
\fm(V)=\pab
\end{equation}
\begin{equation}\label{eqn:reverse}
\fr(V)=(b_\ell,a_\ell,\ldots, b_1,a_1)
\end{equation}
Then $\fm(\fr(V))=\fm(V)^{-1}$. Let $g V g^{-1}$ denote
$(ga_1g^{-1}, gb_1 g^{-1},\ldots,ga_\ell g^{-1}, gb_\ell g^{-1})$ for $g\in G$.
\end{no}

With the above notation, the representation varieties $\fl{\ell}{i}$
can be written as follows:
\begin{eqnarray*}
\fl{\ell}{0}&=&\{V \in G^{2\ell} \mid \fm(V)=e \} \\
\fl{\ell}{1}&=&\{(V,c)\mid V\in G^{2\ell},\ c\in G,\        \fm(V)=c^2 \} \\
\fl{\ell}{2}&=&\{(V,d,c) \mid V\in G^{2\ell},\ d,c\in G,\ \fm(V)=c d c^{-1} d \}
\end{eqnarray*}

\subsection{Symmetric representation varieties of flat connections} \label{sec:Zfl}
Let $\Si$ be a closed, compact, connected, nonorientable surface, and let
$\pi:\tSi\to \Si$ be the orientable double cover. The goal of this and the
next subsection is to relate the representation varieties of $\Si$ to
those of $\tSi$.

Let $\Si=\Si^\ell_i$, where $i=1,2$.
Then $\Si$ is homeomorphic to the connected sum of $2\ell+i$ copies of $\RP^2$,
and its orientable double cover  $\tSi$ is $\Si^{2\ell+i-1}_0$, a Riemann surface
of genus $2\ell+i-1$.

In the rest of this subsection, we follow \cite[Section 5]{Ho} closely.
Define
\begin{eqnarray*}
\Zfl{1}&=&\{(V,c,\bV,\bc)\mid V,\bV\in G^{2\ell},\ c,\bc\in G,\
\fm(V)=c\bc, \fm(\bV)=\bc c\}\\
\Zfl{2}&=&\{(V,d,c,\bV,\bd,\bc)\mid V,\bV\in G^{2\ell},\ d,c,\bd,\bc\in G,\
\fm(V)=c\bd c^{-1}d, \fm(\bV)=\bc d\bc^{-1}\bd \}
\end{eqnarray*}

\begin{lm}\label{thm:ZXfl}
For $i=1,2$,  define $\Phi_G^{\ell,i}: G^{2(2\ell+i)}\to G^{2(2\ell+i-1)}$ by
\begin{eqnarray*}
\Phi_G^{\ell,1}(V,c,\bV,\bc)&=&(V, c\,\fr(\bV)\, c^{-1}),\\
\Phi_G^{\ell,2}(V,d,c,\bV,\bd,\bc)&=&(V, d^{-1}c\, \fr(\bV)\, c^{-1}d, d^{-1},c\bc).
\end{eqnarray*}
where $V,\bV\in G^{2\ell}$, $c,d,\bc,\bd\in G$. Then
$$
\Phi_G^{\ell,i}(\Zfl{i})= \fl{2\ell+i-1}{0}.
$$
\end{lm}
\begin{proof} It is straightforward to check
$\Phi_G^{\ell,i}(\Zfl{i})\subset \fl{2\ell+i-1}{0}$. It remains
to show that $\fl{2\ell+i-1}{0}\subset \Phi_G^{\ell,i}(\Zfl{i})$.
\begin{enumerate}
\item[1.] $\fl{2\ell}{0}\subset \Phi_G^{\ell,1}(\Zfl{1})$:
Given $(V_1,V_2)\in \fl{2\ell}{0}$, where $V_1,V_2\in G^{2\ell}$,
we have $\fm(V_1)\fm(V_2)=e$. Let
$$
\bc=\fm(V_1) =\fm(V_2)^{-1}=\fm(\fr(V_2)).
$$
Then
$$
(V_1,e,\fr(V_2),\bc)\in \Zfl{1},\quad
(V_1,V_2)=\Phi_G^{\ell,1}(V_1,e,\fr(V_2),\bc).
$$

\item[2.] $\fl{2\ell+1}{0}\subset \Phi^{\ell,2}_G(\Zfl{2})$:
Given $(V_1,V_2, a, b)\in \fl{2\ell+1}{0}$,
where $V_1,V_2\in G^{2\ell}$ and $a,b\in G$, we
have $\fm(V_1)\fm(V_2)[a,b]=e$. Let
\[
d=a^{-1},\quad
\bd= a\, \fm(V_1), \quad
 c=a^{-1},\quad
 \bc=ab
\]
Then
$$
(V_1, d,c, \fr(V_2), \bd,\bc)\in \Zfl{2},\quad
(V_1,V_2, a,b)=
\Phi^{\ell,2}_G(V_1, d,c,\fr(V_2), \bd,\bc).
$$
\end{enumerate}
\end{proof}

Let $G^2$ act on $\Zfl{1}$, $\Zfl{2}$ by
\begin{eqnarray*}
(g_1,g_2)\cdot(V,c,\bV,\bc)
&=& (g_1 V g_1^{-1}, g_1 c g_2^{-1}, g_2 \bV g_2^{-1},g_2 \bc g_1^{-1}),\\
(g_1,g_2)\cdot(V,d,c,\bV,\bd,\bc)
&=& (g_1 V g_1^{-1}, g_1 d g_1^{-1}, g_1 c g_2^{-1},
g_2 \bV g_2^{-1},g_2 \bd g_2^{-1} g_2 \bc g_1^{-1}),
\end{eqnarray*}
respectively, where $V,\bV\in G^{2\ell}$ and  $g_1,g_2,c,\bc,d,\bd\in G$.
\begin{lm}\label{thm:ZGGXG}
The surjection $\Phi^{\ell,i}_G: \Zfl{i}\to \fl{2\ell+i-1}{0}$
induces homeomorphisms
\begin{equation}\label{eqn:ZXG}
\Zfl{i}/ G^2 \cong \fl{2\ell+i-1}{0}/G \cong \Hom(\pi_1(\Si^{2\ell+i-1}_0),G)/G
\end{equation}
and a homotopy equivalence
\begin{equation}\label{eqn:ZXhG}
\Zfl{i}^{h G^2} \sim \fl{2\ell+i-1}{0}^{h G}
\end{equation}
where $X^{h G}$ denote the homotopic orbit space $EG\times_G X$.
\end{lm}
\begin{proof} The case $i=1$ of \eqref{eqn:ZXG} was proved in \cite{Ho}; the case $i=2$ is similar.

To see \eqref{eqn:ZXhG}, let $G_1=G\times \{1\}\subset G\times G$ and
$G_2=\{1\}\times G\subset G\times G$, so that $G_1\cong G_2\cong G$,
and the $G_1$-action and $G_2$-action on
$\Zfl{i}$ commute. Note that $G_2$ is a closed normal subgroup of $G\times G$ and that
$G_2$ acts on $\Zfl{i}$ freely, so the natural
projection $\Zfl{i}\to \Zfl{i}/G_2$ induces a homotopy equivalence
$$
\Zfl{i}^{h (G\times G)} \sim (\Zfl{i}/G_2)^{h G_1}.
$$
It is straighforward to check that the surjection
$\Phi^{\ell,i}_G:\Zfl{i}\to \fl{2\ell+i-1}{0}$ descends to
a homeomorphism $\bar{\Phi}^{\ell,i}_G:\Zfl{i}/G_2 \to \fl{2\ell+i-1}{0}$. Moreover,
$\bar{\Phi}^{\ell,i}_G$ is $G$-equivariant with respect to
the $G_1$-action on $\Zfl{i}/G_2$ and the $G$-action on
$\fl{2\ell+i-1}{0}$, so $\bar{\Phi}^{\ell,i}_G$ induces a
homotopy equivalence
$$
(\Zfl{i}/G_2)^{h G_1} \sim \fl{2\ell+i-1}{0}^{h G}
$$
\end{proof}

By \cite[Theorem 3.3]{HL3}, when $\ell>0$ there is a bijection
$$
\pi_0(\Hom(\pi_1(\Si^\ell_0),G)/G)\longrightarrow \pi_1(G_{ss}),
$$
where $G_{ss}=[G,G]$ is the maximal connected semisimple subgroup of $G$.
Since $G^2$ and $G$ are connected, we conclude that
\begin{cor}\label{thm:Zpi}
Suppose that $i=1,2$ and $\ell\geq 0$. Then when $(\ell,i)\neq (0,1)$ there
is a bijection
$$
\pi_0(\Zfl{i}) \longrightarrow \pi_1(G_{ss}).
$$
\end{cor}

\subsection{Involution on symmetric representation varieties of flat connections}
\label{sec:Zinvolution}
In this section, $i=1,2$.

Let $\tau:\Zfl{i}\to \Zfl{i}$ be
the involution defined in \cite{Ho}:
$$
\tau(V,v,\bV,\bv)=
(\bV,\bv,V,v)
$$
where $V,\bV \in G^{2\ell}$ and $v,\bv\in G^i$.
There is an injection $I: \fl{\ell}{i}\to \Zfl{i}$ given by
$(V,v)\mapsto (V,v,V,v)$ such that
$$
I(\fl{\ell}{i})=\Zfl{i}^\tau
$$
where $\Zfl{i}^\tau$ is the fixed locus of the involution $\tau$.
We will show that $\Zfl{i}^\tau$ corresponds to topologically
trivial flat $G$-bundles over the Riemann surface $\Si^{2\ell+i-1}_0$.
To do so, we first recall the definition of the obstruction map,
which detects the topological type of a flat $G$-bundle.

Let $H$ be the connected component of the identity
of the center of $G$, and let $G_{ss}=[G,G]$ be the commutator group.
Then $H$ is a compact torus, and $G_{ss}$ is the maximal
connected semisimple subgroup of $G$. Let $\rho_{ss}:\tG_{ss}\to G_{ss}$  and 
$\rho:\tG\to G$ be
universal coverings. Then $\tG=\mathfrak{h}\times
\tG_{ss}$ where $\mathfrak{h}=Lie(H)$.
Define  $\rho^{2\ell}:\tG^{2\ell}\to G^{2\ell}$ by 
$$
\rho^{2\ell}(\tab)\mapsto
(\rho(\ta_1), \rho(\tb_1),\ldots, \rho(\ta_\ell), \rho(\tb_\ell)).
$$
With the above notation, the obstruction map
$o:\fl{\ell}{0}\to \Ker(\rho_{ss})\cong \pi_1(G_{ss})$ is
defined as follows:
given $V\in G^{2\ell}$, pick $\tV\in \tG^{2\ell}$ such 
that $\rho^{2\ell}(\tV)=V$, and define $o(V)=\fm(\tV)$.
Then $o(V)\in \Ker\rho \cap G_{ss}=\Ker\rho_{ss}$,
and the definition is independent of choice of
$\tV$. The flat $G$-bundle associated to $V$ is
topologically trivial if and only
if $o(V)=\te$, where $\te$ is the identity element of $\tG$.

Let $o:\fl{2\ell+i-1}{0} \to \Ker(\rho_{ss})$ be the obstruction map, and let
$$
o'=o\circ \Phi_G^{\ell,i}: \Zfl{i} \to \pi_1(G_{ss}).
$$
Let $e$ be the identitiy element of $G$.
\begin{lm}
$o'(\tau(y))=o'(y)^{-1}$ for $y\in \Zfl{i}$.
\end{lm}
\begin{proof} 
We will prove the case $i=1$. The case $i=2$ is similar.

Given $y=(V,c,\bV, \bc)\in \Zfl{1}$, where $V,\bV\in G^{2\ell}$
and $c,\bc\in G$, pick  $V', \bV'\in \tG^{2\ell}$ and
$c', \bc'\in \tG$ such that
$\rho^{2\ell}(V')=V$, $\rho^{2\ell}(\bV')=\bV$,
$\rho(c')=c$, and $\rho(\bc')=\bc$. Then 
$$
\rho(\fm(V')(c'\bc')^{-1}) =\fm(V)(c\bc)^{-1}=e,
\quad \rho(\fm(\bV')(\bc'c')^{-1})=\fm(\bV)(\bc c)^{-1}=e.
$$
Let $k=\fm(V')(c'\bc')^{-1}$ and $\bar{k}=\fm(\bV')(\bc' c')^{-1}$.
Then $k, \bar{k}\in \Ker\rho \subset Z(\tG)$.
We have
\begin{eqnarray*}
o'(y) &=& o'(V,c,\bV,\bc)=o(V, c\fr(\bV) c^{-1})= \fm(V')\fm(c'\fr(\bV') (c')^{-1})\\
&=&\fm(V') c' \fm(\bV')^{-1} (c')^{-1}
= (k c'\bc') c' (\bar{k}\bc' c')^{-1} (c')^{-1}= k\bar{k}^{-1},\\
o'(\tau(y))&=&o'(\bV,\bc,V,c)= \bar{k} k^{-1}.
\end{eqnarray*}
So $o'(\tau(y))=o'(y)^{-1}$.
\end{proof}

\begin{lm}\label{thm:pullback-trivial}
 $o' \circ I(x)= \te$ for all $x\in \fl{\ell}{i}$. 
\end{lm}
\begin{proof} We will prove the case $i=1$. The
case $i=2$ is similar.

Given $(V,c)\in \fl{\ell}{1}$, where $V\in G^{2\ell}$ and
$c\in G$, pick $\tV\in \tG^{2\ell}$ and $\tc\in \tG$
such that $\rho^{2\ell}(\tV)=V$ and $\rho(\tc)=c$.
Then $\rho(\fm(\tV)\tc^{-2}) =\fm(V)c^{-2}=e$. Let $k=\fm(\tV) \tc^{-2}\in 
\Ker\rho\subset Z(\tG)$. Then
\begin{eqnarray*}
&&o'\circ I(V,c)= o'(V,c,V,c)
=o(V,c \fr(V) c^{-1})
=\fm(\tV)\fm(\tc \fr(\tV) \tc^{-1})\\
&& = \fm(\tV) \tc \fm(\tV)^{-1} \tc^{-1} = k\tc^2 \cdot \tc (k\tc^2)^{-1}\tc^{-1}=\te.
\end{eqnarray*}
\end{proof}

By \cite[Theorem 5.2]{HL3}, any topological principal $G$-bundle on a closed, connected,
nonorientable surface admits a flat connection. 
By Lemma \ref{thm:pullback-trivial},
the pullback of a flat $G$-bundle over $\Si^{\ell}_i$
under the orientable double cover $\Si^{2\ell+i-1}_{0}\to 
\Si^{\ell}_{i}$ is
 a topologically trivial flat $G$-bundle over $\Si^{2\ell+i-1}_0$.
We conclude that:
\begin{pro}\label{thm:top_trivial}
Let $G$ be a compact, connected Lie group.
Let $\Si$ be a closed, connected, nonorientable surface, and
let $\pi:\tSi\to \Si$ be the orientable double cover.
Then the pullback
$\pi^*P$ of any topological principal $G$-bundle $P\to \Si$ is
topologically trivial.
\end{pro}

\section{Yang-Mills Functional and Yang-Mills Connections}\label{sec:functional}

In this section, we will define Yang-Mills functional and Yang-Mills
connections on nonorientable closed surfaces.

\subsection{Yang-Mills functional and Yang-Mills connections on orientable surfaces}
We first recall the Yang-Mills functional and Yang-Mills connections on orientable closed
surfaces, following \cite{ym}.

Let $G$ be a compact connected Lie group.
Let $\Si$ be a Riemann surface. There is a unique K\"{a}hler
metric $h$ such that the scalar curvature is a constant and
the K\"{a}hler form $\omega$ is the unique harmonic 2-form on
$\Si$ such that $\int_\Si\omega=1$. We call it the {\em canonical metric}
of the Riemann surface.

 Let $\cA(P)$ denote the space of $C^\infty$ connections on $P$.
Then $\cA(P)$ is an affine space whose  associated (real) vector space is
$\Omega^1(\Si,\ad(P))$.   The {\em Yang-Mills functional} $L:\cA(P)\to \bR$
is defined by
\begin{equation}\label{eqn:L}
L(A)=\int_\Si \Tr (F(A)\wedge * F(A))
\end{equation}
where $F(A)$ is the curvature form of $A$.

Let $A_t=A+t\eta$ be a line of connections, where
$\eta\in \Omega^1(\Si,\ad(P))$. Then
$$
F(A_t)=F(A)+t d_A \eta +\frac{1}{2} t^2[\eta,\eta],
$$
so
\begin{eqnarray*}
L(A_t) &=& L(A)+ 2t \int_\Si \Tr ( d_A\eta \wedge * F(A)) + O(t^2)\\
&=& L(A)+2t \int_\Si \Tr( \eta \wedge * d_A * F(A) )+ O(t^2)
\end{eqnarray*}
So $A$ is a critical point of $L$ iff it satisfies the {\em Yang-Mills
equation}:
\begin{equation}\label{eqn:YM}
d^*_A F(A)=* d * F(A)=0.
\end{equation}
We call critical points {\em Yang-Mills connections} on $P$.
Note that flat connections are Yang-Mills connections.

\subsection{Involution on the principal bundle}\label{sec:tau_P}

Let $\Si$ be a connected, nonorientable, closed surface. Then $\Si$ is diffeomorphic to
the connected sum of $m>0$ copies of $\RP^2$'s.
Let $\pi: \tSi \to \Si$ be the orientable double cover, and
let $\tau:\tSi \to \tSi$ be the deck transformation. Then
$\tSi$ is a Riemann surface of genus $m-1$, and
$\tau$ is an anti-holomorphic, anti-symplectic involution with no fixed point.

Let $P\to \Si$ be a principal $G$-bundle. Let $\tP=\pi^* P$ be the pullback
principal $G$-bundle on $\tSi$. By Proposition \ref{thm:top_trivial}, $\tP$ is topologically trivial.
There is an involution $\tilde{\tau}: \tP\to \tP$ which is
$G$-equivariant and covers $\tau:\tSi\to \tSi$.

More explicitly, fix a trivialization $\tP\cong \tSi\times G$.
The right $G$-action on $\tP$ is given by
$$
(x,h)\cdot g = (x,h\cdot g)
$$
where $g\in G, (x,h)\in \tSi\times G$.
It is straightforward to check that
$$
(x,h)\cdot (g_1g_2)= ((x,h)\cdot g_1)\cdot g_2.
$$
The involution $\tilde{\tau}$ is $G$-equivariant with respect to the above
$G$-action:
$$
\tilde{\tau}(x,h)\cdot g = \tilde{\tau}((x,h)\cdot g)
$$
for $(x,h)\in\tP$, $g\in G$. Let $s:\tSi\to G$ be defined by
$$
\tilde{\tau}(x,e)=(\tau(x), s(x))
$$
where $e\in G$ is the identity element. By the $G$-equivariance,
\begin{equation}\label{eqn:tau_s}
\tilde{\tau}(x,h)=(\tau(x), s(x)h).
\end{equation}
We have $\tilde{\tau}\circ \tilde{\tau}=\mathrm{id_{\tP}}$, so
\begin{equation}\label{eqn:s}
s(\tau(x))=s(x)^{-1}.
\end{equation}

Conversely, given any continuous map $s:\tSi\to G$ such that
\eqref{eqn:s} holds, we define $\tilde{\tau}_s:\tP\to \tP$ by
\eqref{eqn:tau_s}. Then $\tilde{\tau}_s$ is a $G$-equivariant involution
on $\tP$ which covers the involution $\tau$ on $\tSi$, so
$P_s=\tP/\tilde{\tau}_s$ is a principal $G$-bundle over $\Si$.

In particular, we can take $s$ to be a constant map:
$s(x)\equiv \ep$, where $\ep\in G$, $\ep^2=e$.
The involution $\tilde{\tau}_\ep \equiv \tilde{\tau}_s$
on $\tP\cong \tSi\times G$ is given by
$$
(x,h)\mapsto (\tau(x),\ep h).
$$
The zero connection on $\tSi\times G$ descends to
a flat connection $A_\ep$ on $P_\ep =\tP/\tilde{\tau}_\ep$
which corresponds to
$$
(e,\ldots,e, \epsilon)\in \fl{\ell}{i}\subset G^{2\ell+i}.
$$

The topological type of $P_\ep=\tP/\tilde{\tau}_\ep$ can
be determined by the following way (see \cite{HL3}).
We use the notation in Section \ref{sec:Zinvolution}.
Choose $\tilde{\epsilon}\in \rho^{-1}(\epsilon)$,
where $\rho:\tG\to G$ is the universal covering. The
{\em obstruction class}
$$
[\tilde{\epsilon}^2] \in \Ker\rho/2\Ker\rho
\cong \pi_1(G)/2\pi_1(G)\cong H^2(\Si;\pi_1(G))
$$
is independent of the choice of $\tilde{\epsilon}$ and determines
the topological type of $P_{\epsilon}$. Recall that
$\Ker\rho\cong \pi_1(G)$ is abelian, and
$$
\Prin_G(\Si)\cong H^2(\Si;\pi_1(G)).
$$

Conversely, a principal $G$-bundle over $\Si$ of any topological
type arises this way. Recall that
$$
\Ker\rho\subset Z(\tilde{G})\subset \mathfrak{h}\times \tilde{T}_{ss}
$$
where $\tilde{T}_{ss}$ is some maximal torus of $\tilde{G}_{ss}$.
Given $k\in \Ker\rho/2\Ker\rho$ represented by
$a\in \Ker\rho$, choose $\tilde{\epsilon}\in
\mathfrak{h}\times \tilde{T}_{ss}$ such that $\tilde{\epsilon}^2 =a$.
Let $\epsilon=\rho(\tilde{\epsilon})\in G$.
Then $\epsilon^2=e$, and  $\epsilon$ defines a principal $G$-bundle $P_\epsilon\to \Si$
with obstruction class $k$.

\begin{ex}
$G=U(n)$. Let $\epsilon\in U(n)$  such that $\epsilon^2= I_n$ (in particular,
$\det(\epsilon)=\pm 1$). Then $c_1(P_\epsilon)=c_1(\det(P_\epsilon))$, where $\det(P_\epsilon)$
is the $U(1)$-bundle on $\Si$ which is the quotient of $\tSi\times U(1)$
by the involution
$$
(x,h)\mapsto (\tau(x), \det(\epsilon) h)
$$
for $x\in\tSi$, $h\in U(1)$. So $P_\epsilon \cong P_\pm$ if $\det(\epsilon)=\pm 1$, where
$c_1(P_+)=0$ and $c_1(P_-)=1$ in $H^2(\Si;\bZ)\cong \bZ/2\bZ$.
\end{ex}

\subsection{Involution on the adjoint bundle}\label{sec:tau_ad}
Let
$\fg$ denote the Lie algebra of $G$.
Let $\tP=\tSi\times G$ be the trivial principal
$G$-bundle as above. Let
$\ad(\tP)=\tP\times_{G}\fg$,
where $G$ acts on $\tP\times \fg$ by
$$
g\cdot (x,h,X) = (x, hg, \Ad(g^{-1})(X))
$$
for $g\in G$, $(x,h)\in \tSi\times G=\tP$, $X\in \fg$.
Then $\ad(\tP)\cong \tSi\times \fg$, and the natural projection
$\tP\times \fg\to \tSi\times \fg$ is given by
$$
(x,h,X)\mapsto (x,\Ad(h)(X)).
$$

Let $s: \tSi\to G$ be a smooth map such that
$s(\tau(x))=s(x)^{-1}$, as in \eqref{eqn:s};
define $\tilde{\tau}_s:\tP\to \tP$ by
$\tilde{\tau}_s(x,h)=(\tau(x), s(x)h)$, as in \eqref{eqn:tau_s}.
The involution $\tilde{\tau}_s$ on $\tP$ induces an involution
on $\ad(\tP)$:
$$
(x,X)\mapsto (\tau(x), \Ad(s(x))(X)).
$$
We use the same notation $\tilde{\tau}_s$ to denote it.
We have
$$
\ad(\tP)/\tilde{\tau}_s\cong P_s \times_G\fg.
$$

\subsection{Involution on the space of connections}\label{sec:tau_A}
$\tilde{\tau}_s:\ad(\tP)\to \ad(\tP)$ induces an involution
$$
\tilde{\tau}_s^*:\cA(\tP)\to \cA(\tP),
$$
where $\cA(\tP)\cong \Omega^1(\tSi;\fg)$.

More explicitly,
given $\theta\in \Omega^1(\tSi)$ and $X\in \Omega^0(\tSi;\fg)$,
$$
\tilde{\tau}_s ^*(X\otimes \theta) =
\Ad(s)(\tau^* X)\otimes\tau^* \theta.
$$
Similarly, $\tilde{\tau}_s:\ad(\tP)\to \ad(\tP)$ induces an involution
on $\Omega^2(\tSi;\fg)$.
The curvature form $F(A)$ can be viewed as an element in $\Omega^2(\tSi;\fg)$:
$$
F(A)=X\otimes \omega
$$
where $X:\tSi\to \fg$ and $\omega$ is the volume form of $\tSi$.
We have
$$
F(\tilde{\tau}_s^* A)= \tilde{\tau}_s^*(F(A))
=\tilde{\tau}_s^*(X\otimes \omega) =\Ad(s)(\tau^* X)\otimes\tau^*\omega
=-\Ad(s)(\tau^* X)\otimes\omega
$$
where we have used the fact that $\tau$ is anti-symplectic.

Recall that $\cA(\tP)$ is a K\"{a}hler manifold: the complex structure
is given by $\alpha\mapsto * \alpha$, and the symplectic form $\Omega$ is given
by
$$
\Omega(\alpha,\beta)=\int_{\tSi} \Tr(\alpha\wedge \beta).
$$
The involution $\tilde{\tau}_s^*:\cA(\tP)\to \cA(\tP)$ is
anti-holomorphic and anti-symplectic.
The fixed locus $\cA(\tP)^{\tau_s}$ can be identified with
$\cA(P_s)$, the space of $G$-connections on $P_s$. $\cA(\tP)^{\tau_s}$ is a
totally geodesic, totally real, Lagrangian submanifold of $\cA(\tP)$.

\subsection{Yang-Mills functional and Yang-Mills connections on nonorientable surfaces}

Let $(M,g)$ be a Riemannian manifold with an isometric involution
$\tau:M\to M$. It is straightforward to check the following statements.

\begin{lm}\label{thm:grad}
Let $f: M\to \bR$ be a smooth function
such that $f\circ\tau=f$.
\begin{enumerate}
\item Let $N$ be the set of critical points of $f$. Then $\tau(N)=N$.
\item  Let $X$ be the gradient vector field of $f$.
Then
\begin{enumerate}
\item For any $p\in M$, we have $\tau_*(X(p))= X(\tau(p))$.
\item If $\gamma:I \to M$ is an integral curve
of $X$, where $I$ is an open subset of $\bR$, so is
$\tau\circ \gamma:I\to M$.
\end{enumerate}
\end{enumerate}
\end{lm}

Let $M^\tau$ be the fixed locus of $\tau$. Suppose that
$$
M^\tau=\bigcup_{i\in I} M^\tau_i
$$
is a union of connected components, where each $M^\tau_i$ is
a submanifold of $M$. Then each $M^\tau_i$ is a totally
geodesic submanifold of $M$ because $\tau$ is an isometry.
It is straightforward to check the following statements.
\begin{lm}
Let $f:M\to \bR $ be a smooth function such that
$f\circ \tau = f$, and let $f^\tau: M^\tau \to \bR$ be the
restriction of $f$.Then
\begin{enumerate}
\item $X(p)\in T_p(M^\tau)$ for any $p\in M^\tau$, and
          $X|_{M^\tau }$ is the gradient vector field of $f^\tau$.
\item The set of critical points of $f^\tau: M^\tau\to \bR$
is $N^\tau=N\cap M^\tau$, where $N$ is the set of
critical points of $f:M\to \bR$.
\end{enumerate}
\end{lm}

In our case,  $M=\cA(\tP)$ and $f$ is the Yang-Mills functional $L$.
We define the Yang-Mills functional on $\cA(\tP)^{\tau_s}\cong\cA(P_s)$ to be
$L^{\tau_s}:\cA(\tP)^{\tau_s}\to \bR$. We call the critical points of $L^{\tau_s}$ {\em
Yang-Mills connections} on $P$. By Lemma \ref{thm:grad},
$A$ is a Yang-Mills connection on $P$ if and only if
$\pi^*A$ is a Yang-Mills connection on $\tP$.

It is worth mentioning that our definition of Yang-Mills connections on non-orientable surfaces
is different from the one introduced by S. Wang in \cite{Wa}.

\section{Yang-Mills Connections and Representations}\label{sec:YMrep}
In this section,
we introduce the {\em super central extension} of the fundamental
group of a surface, and
establish a precise correspondence between Yang-Mills connections and
representations of super central extension. We also
introduce representation  varieties for Yang-Mills connections
for orientable and nonorientable surfaces, 
and extended moduli spaces for nonorientable surfaces.

\subsection{Super central extension of the fundamental group}
To relate Yang-Mills connections to representations, we need to introduce certain
extension of the fundamental group of the surface.

Let $\Si$ be a closed, compact, connected surface.
Given $a\in \pi_1(\Si)$, let $\deg(a)=w_1(T_\Si)[a]\in \bZ/2\bZ$,
where $[a]$ is the image of $a$ under the group homomorphism
$\pi_1(\Si)\to H_1(\Si;\bZ)=\pi_1(\Si)/[\pi_1(\Si),\pi_1(\Si)]$, and
$w_1(T_\Si)\in H^1(\Si;\bZ/2\bZ)$ is the first Stiefel-Whitney class of the tangent bundle
of $\Si$. More geometrically, if $\gamma:S^1\to \Si$ is a loop representing $a\in \pi_1(\Si)$,
then $\deg(a)=0\in \bZ/2\bZ$ if the rank 2 real vector bundle $\gamma^* T_\Si$ over $S^1$ is
orientable (or equivalently, topologically trivial); $\deg(a)=1 \in \bZ/2\bZ$
if $\gamma^* T_\Si$ is non-orientable (or equivalently, topologically non-trivial).
The group homomorphism $\deg:\pi_1(\Si)\to \bZ/2\bZ$ is trivial if and only if $\Si$ is orientable.

We are now ready to define the {\em super central extension} $\Gamma_\bR (\Si)$ of $\pi_1(\Si)$.
It fits in a short exact sequence of groups:
$$
1 \to \bR \stackrel{\alpha}{\to} \Gamma_\bR(\Si)
\stackrel{\beta}{\to} \pi_1(\Si) \to 1.
$$
Given $r\in \bR$, let $J_r=\alpha(r)$, so that
$J_{r_1+r_2}=J_{r_1} J_{r_2}$. Given $a\in \Gamma_\bR(\Si)$, we have
$$
a J_r a^{-1} = \left\{
\begin{array}{ll}
J_r \\ J_{-r}=J_r^{-1}
\end{array}
\right.
\ \textup{ if }\  \deg(\beta(a))=
 \left\{
\begin{array}{ll}
0\\ 1
\end{array}
\right.
\in \bZ/2\bZ.
$$
This defines $\Gamma_\bR(\Si)$ up to group isomorphism.
We will give a more explicit description later.

 When $\Si$ is orientable, $\Gamma_\bR(\Si)$ is the  central
extension of $\pi_1(\Si)$ defined in \cite{ym}.

\subsection{Representation varieties for orientable surfaces}
Recall that any closed, compact, connected surface
is diffeomorphic to $\Si^\ell_0$, a Riemann surface of genus $\ell$, for some
nonnegative integer $\ell$. $\Gamma_\bR(\Si^\ell_0)$ is generated by
$$
A_1,B_1,\ldots, A_\ell, B_\ell, J_r
$$
where $r\in\bR$, with relations
\begin{enumerate}
\item $J_{r_1}J_{r_2}= J_{r_1+r_2}$
\item $[A_i ,J_r]=[B_i,J_r]=1,\quad  i=1,\ldots,\ell,\quad r\in\bR$.
\item $\prod_{i=1}^\ell [A_i,B_i] = J_1$.
\end{enumerate}

Let $\rho: \Gamma_\bR(\Si^\ell_0)\to G$ be a group homomorphism.
From the relation (1) we must have $\rho(J_r)=\exp(rX)$
for some $X\in \fg$, where $\fg$ is the Lie algebra of $G$. From
the relation (2) we must have $\rho(A_i), \rho(B_i)\in G_X$,
where $G_X$ is the stabilizer of $X$ of the adjoint action of $G$ on $\fg$.
Combined with the relation (3), $\Hom(\Gamma_\bR(\Si^\ell_0),G)$ can be identified
with
\[
\ym{\ell}{0}=\{ (V,X)\in G^{2\ell}\times \fg\mid
V\in (G_X)^{2\ell} , \fm(V)=\exp(X) \}
\]
where $\fm(\ab)=\pab$ was defined in Section \ref{sec:Xfl}.

Let $\cN(P)\subset \cA(P)$ be the space of Yang-Mills connections on
$P$. By Theorem $6.16$ in \cite{ym}, $\cN(P)$ is nonempty for any underlying
principal $G$-bundle $P$. Let $\cN_0(P)\subset \cA(P)$ be the space of flat connections
on $P$, as in Section \ref{sec:flat}.
The natural inclusion $\cN_0(P)\subset \cN(P)$ induces
an inclusion $\fl{\ell}{0}\hra \ym{\ell}{0}$, $V\mapsto (V,0)$.

Let $\cG(P)$ be the gauge group, and let $\cG_0(P)$ be the based gauge group, as before.
\begin{thm}[{\cite[Theorem 6.7]{ym}}]\label{thm:YMzero}
There is a bijective correspondence between conjugacy classes of 
homomorphisms $\Gamma_\bR(\Si)$
and gauge equivalence classes of Yang-Mills $G$-connections over $\Si$. In other words, we have
\begin{eqnarray*}
\bigcup_{P\in \Prin_G(\Si^\ell_0)} \cN(P)/\cG_0(P)
&\cong& \Hom(\Gamma_\bR(\Si^\ell_0),G)\cong \ym{\ell}{0},\\
\bigcup_{P\in \Prin_G(\Si^\ell_0)} \cN(P)/\cG(P)
&\cong&\Hom(\Gamma_\bR(\Si^\ell_0),G)/G \cong \ym{\ell}{0}/G,
\end{eqnarray*}
where $g\in G$ acts on $G^{2\ell}\times \fg$ by
$\ g\cdot (V,X)= (gVg^{-1}, \Ad(g)(X) )$.
\end{thm}

\subsection{Holonomy on the double cover}\label{sec:holonomy}
Let $\Si$ be a closed, compact, connected, nonorientable
surface, and let $\pi:\tSi\to \Si$ be the orientable
double cover. Then $\tSi$ is a closed, compact, connected
orientable surface. Let $\tau:\tSi\to\tSi$ be the deck
transformation which is an orientation reversing involution.

Let $A\in\cA(\tP)$ be a Yang-Mills connection.
Recall that $A\in \cA(\tP)$ is a Yang-Mills connection if and only if there exists
$u:\tSi\to G$ such that
$$
F(A) =\Ad(u)(X)\otimes \omega
$$
where $X$ is a {\em constant} vector in $\fg$, or equivalently, if there
exists $A'\in \cA(\tP)$ and $u\in\cG(\tP)$ such that
$A=u\cdot A'$ and $F(A')= X\otimes \omega$.

We fix
a trivialization $\tP\to \tSi\times G$ such that
$$
F(A)=X\otimes \omega
$$
where $X\in\fg$ and $\omega$ is the volume form.
Using this trivialization, we may define the holonomy along
a {\em path} (the holonomies along based loops are defined without using
the trivialization of $\tP$). Given a path $\gamma:[0,1]\to \tSi$,
let $\tilde{\gamma}:[0,1]\to \tP$ be the horizontal lifting of $\gamma$
(with respect to the connection $A$) with $\tilde{\gamma}(0)=(\gamma(0),e)$,
where $e$ is the identity element. Then $\tilde{\gamma}(1)=(\gamma(1),g^{-1})$
for some $g\in G_X$, where $G_X$ is the stabilizer
of $X\in \fg$ of the adjoint action of $G$ on $\fg$.
We call $g\in G$ the holonomy along $\gamma$.

Let $\tilde{\gamma}':[0,1]\to \tP$ be another horizontal lifting
of $\gamma$ with $\tilde{\gamma}'(0)=(\gamma(0),h)$,
where $h\in G$. By $G$-invariance of the connection, we have
$\tilde{\gamma}' =\tilde{\gamma}\cdot h$, so
$$\tilde{\gamma}'(1)=(\gamma(1), g^{-1}h)=(\gamma(1), h(h^{-1}g h)^{-1} ).
$$
To summarize, if we change the trivialization by a constant gauge
transformation $h$, the curvature form changes from
$F(A)=X\otimes \omega\in\Omega^2(\tSi;\fg)$ to
$$
F(A)= \Ad(h^{-1})X\otimes \omega\in \Omega^2(\tSi;\fg)
$$
and the holonomy along $\gamma$ changes from $g\in G_X$ to
$h^{-1}gh \in G_{\Ad(h^{-1})X}$.

Recall that $\Si$ is diffeomorphic to $\Sigma^\ell_i$ for some $\ell\geq 0$ and
$i=1,2$, where $\Si^\ell_1$ is the connected sum of a Riemann surface of genus
$\ell$ and the real projective plane, and $\Si^\ell_2$ is the connected sum
of a Riemann surface of genus $\ell$ and a Klein bottle. We will discuss the
case $\Si^\ell_1$ in detail. The case $\Si^\ell_2$ is similar.

Suppose that $s:\tilde{\Sigma}\to G$ satisfies \eqref{eqn:s}, so it
defines an involution $\tilde{\tau}_s:\tP\to\tP$.
Now look at Figure 1.

\begin{figure}[!htbp]
        \begin{center}
                \psfrag{a2}[c][c][1][0]{$\bar{A}$}
                \psfrag{b2}[c][c][1][0]{$\bar{B}$}
                \psfrag{c2}[c][c][1][0]{$\bar{C}$}
                \psfrag{c1}[c][c][1][0]{$C$}
                \psfrag{a1}[c][c][1][0]{$A$}
                \psfrag{b1}[c][c][1][0]{$B$}
                \psfrag{p1}[c][c][1][0]{$p_{+}$}
                \psfrag{p2}[c][c][1][0]{$p_{-}$}
                \includegraphics[scale = 0.7]
                {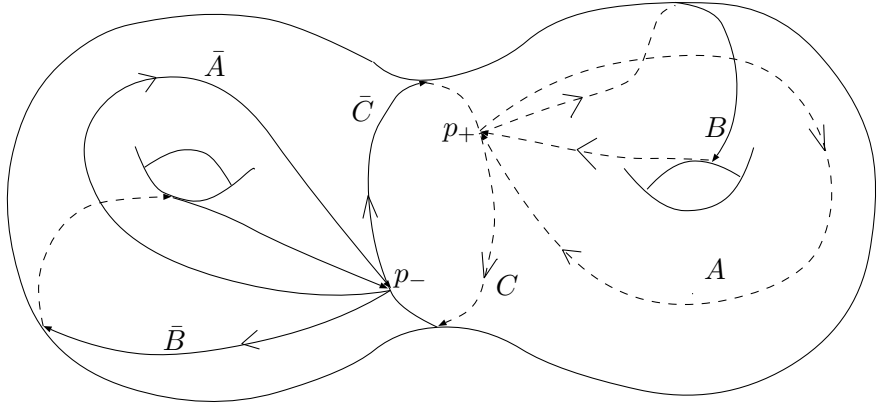}
        \end{center}
        \caption{Holonomy on the double cover $\Si^2_0$}
\end{figure}

$A_i,B_i$ are loops passing
through $p_+\in\tSi$, $\bar{A}_i,\bar{B}_i$ are loops passing through
$p_-=\tau(p_+)$, $C$ is a path from $p_+$ to $p_-$, and $\bar{C}$ is a path
from $p_-$ to $p_+$. The holonomies along
$A_i,B_i,\bar{A}_i,\bar{B}_i$ depend on the connection,
not on the trivialization. The holonomies along
$C$ and $\bar{C}$ depend on the connection and the trivialization. We choose
the trivialization as follows. Let the trivialization of $\tP$ at
$p_+$ and $p_-$ be given by $h\mapsto h$ and $h\mapsto s(p_+) h$,
respectively. We define $c$ and $\bc$ as follows.
Let $\gamma:[0,1]\to \tP$ be the horizontal lifting of $C$ such that
$\gamma(0)=(p_+,e)$. Then
$\gamma(1)=(p_-, s(p_+)c^{-1})$.
Let $\bar{\gamma}:[0,1]\to \tP$ be the horizontal lifting of $\bar{C}$
such that $\bar{\gamma}(0)=(p_-, s(p_+))$. Then $\bar{\gamma}(1)=(p_+,\bc^{-1})$.
Let $\gamma'(t)=\bar{\gamma}(t)\cdot c^{-1}$. Then
$\gamma'$ is also a horizontal lifting of $\bar{C}$,
$\gamma'(0)= (p_-, s(p_+)c^{-1})$, $\gamma'(1)=(p_+, \bc^{-1} c^{-1})=(p_+,(c\bc)^{-1})$.
So $\gamma \cup \gamma'$ is a horizontal lifting of $C\bar{C}$,
and the holonomy along $C \bar{C}$ is $c\bc$. Denote the holonomies
along $A_i,B_i,\bar{A}_i,\bar{B}_i$
by $a_i,b_i,\bar{a}_i,\bar{b}_i$ respectively.

We cut $\tSi$ into two discs $D_+$ and $D_-=\tau(D_+)$. 
The (oriented) boundaries of
of $D_+$ and $D_-$ are
$$
\partial D_+= \prod_{i=1}^\ell [\bar{A}_i,\bar{B}_i]C^{-1} \bar{C}^{-1},\quad\quad
\partial D_-=\prod_{i=1}^\ell [A_i,B_i] \bar{C}^{-1} C^{-1}.
$$
Recall that
$$
\int_{\tSi}\omega =1, \quad
\int_{\tSi}\omega = \int_{D_+} \omega +\int_{-D_-}\omega, \quad
\tau^* \omega=-\omega
$$
where $-D_-$ is $D_-$ with the reversed orientation. We conclude that
$$
\int_{D_\pm}\omega =\pm \frac{1}{2}.
$$

Let $z=s(p_+)\in G$.
From the above discussion, we have
\begin{eqnarray*}
\fm(V) \bc^{-1}c^{-1}&=&
\exp\Bigl(-\int_{D_-}X\otimes \omega\Bigr)=\exp(X/2)\\
\fm(\bV) c^{-1} \bc^{-1}
&=&\exp\Bigl(-\int_{D_+}\Ad(z^{-1})X\otimes \omega\Bigr)=\exp(-\Ad(z^{-1})X/2)
\end{eqnarray*}

Moreover,
$$
V \in (G_X)^{2\ell}, cz ^{-1} \in G_X,\quad
\bV\in (G_{\Ad(z^{-1})(X)})^{2\ell},\quad
\bc z  \in G_{\Ad(z^{-1})(X)}
$$

So we shall define a symmetric representation variety
\begin{eqnarray*}
\Zym{1}^z &=& \Bigl\{(V,c,\bV,\bc,X)\in G^{2(2\ell+1)}\times \fg \mid
V\in (G_X)^{2\ell},\ c z^{-1} \in G_X,\\
&& \quad \bV\in (G_{\Ad(z^{-1}(X)})^{2\ell},\
\bc z  \in G_{\Ad(z^{-1})(X)},\\
&&\fm(V) =\exp(X/2)c\bc,\
\fm(\bV) =\exp(-\Ad(z^{-1})X/2) \bc c \Bigr\}
\end{eqnarray*}

Our next goal is to rewrite $\Zym{1}^z$ without
using $z$. Given
$(V,c,\bV,\bc,X)\in \Zym{1}^z$, we have
$\bc z  \in G_{\Ad(z^{-1})(X)}$, which implies
$\Ad(z^{-1})(X)=\Ad(\bc)(X)$. So
$$
\fm(\bV)=\exp(-\Ad(\bc)(X)/2)\bc c
=\bc \exp(-X/2)c.
$$
We also have
$$
c\bc\in G_X,\quad
\bV\in (G_{\Ad(\bc)(X)})^{2\ell},\quad
\bc c\in G_{\Ad(\bc)(X)},
$$
which imply
$$
c\bV c^{-1}\in (G_X)^{2\ell}.
$$
Define
\begin{eqnarray*}
\Zym{1}&=&\bigl\{(V,c,\bV,\bc,X)\in G^{2(2\ell+1)}\times \fg\mid
V,c\bV c^{-1} \in (G_X)^{2\ell}, \\
&& \quad \fm(V)=\exp(X/2)c\bc, \ \fm(\bV)=\bc \exp(-X/2)c\bigr\}
\end{eqnarray*}
Then
$$
\Zym{1}^z =\{ (V,c,\bV,\bc,X)\in \Zym{1} \mid \Ad(\bc)(X)=\Ad(z^{-1})(X) \},
$$
where $V,\bV\in G^{2\ell}, c,\bc \in G, X\in \fg$, and
$$
\Zym{1}=\bigcup_{z\in G} \Zym{1}^z.
$$

The involution $\tilde{\tau}_s:\cA(\tP)\to \cA(\tP)$ induces a map
$\tau_z: \Zym{1}^z \to \Zym{1}^{z^{-1}}$ given by
$$
(V,c,\bV,\bc,X)\mapsto (\bV,\bc,V,c,-\Ad(z^{-1})X),
$$
or equivalently,
\begin{equation}\label{eqn:inv_1}
(V,c,\bV,\bc,X)\mapsto (\bV,\bc,V,c,-\Ad(\bc)X).
\end{equation}
Note that \eqref{eqn:inv_1} defines an involution $\tau: \Zym{1}\to \Zym{1}$.
Similarly,
the involution $\tilde{\tau}_s:\cA(\tP)\to \cA(\tP)$ induces an involution
$\tau: \Zym{2}\to \Zym{2}$
given by
\begin{equation}\label{eqn:inv_2}
(V,d,c,\bV,\bd,\bc,X)\mapsto (\bV,\bd,\bc,V,d,c,-\Ad(\bc)X).
\end{equation}

\subsection{Symmetric representation varieties}
In this subsection, $i=1,2$. Based on the discussion in Section \ref{sec:holonomy},
we define symmetric representation varieties as follows:
\begin{eqnarray*}
\Zym{1}&=&\bigl\{(V,c,\bV,\bc,X)\in G^{2(2\ell+1)}\times \fg\mid
 V , c\bV c^{-1}\in (G_X)^{2\ell}, \\
&&\quad \fm(V)=\exp(X/2)c\bc, \ \fm(\bV)=\bc \exp(-X/2) c \bigr\}\\
\Zym{2}&=&\bigl\{(V,d,c,\bV,\bd,\bc,X)\in G^{2(2\ell+2)}\times \fg\mid
V , d^{-1}c\,\bV\, c^{-1}d \in (G_X)^{2\ell}, \quad d^{-1},c\bar{c}\in G_X\\
&& \quad \fm(V)=\exp(X/2)c\bd c^{-1}d, \ \fm(\bV)=\bc d\exp(-X/2) \bc^{-1}\bd \bigr\}
\end{eqnarray*}
where $\fm(\ab)=\pab$ was defined in Section \ref{sec:Xfl}.
\begin{lm}\label{thm:ZXym}
For $i=1,2$,  define $\Phi_G^{\ell,i}: G^{2(2\ell+i)}\times \fg
\to G^{2(2\ell+i-1)}\times \fg$ by
\begin{eqnarray*}
\Phi_G^{\ell,1}(V,c,\bV,\bc,X)&=&(V, c\,\fr(\bV)\, c^{-1},X)\\
\Phi_G^{\ell,2}(V,d,c,\bV,\bd,\bc,X)
&=&(V, d^{-1}c\,\fr(\bV)\, c^{-1}d, d^{-1},c\bc,X)
\end{eqnarray*}
where $V,\bV\in G^{2\ell}$, $c,d\in G$, $X\in\fg$, and
$\fr(\ab)=(b_\ell,a_\ell,\ldots,b_1,a_1)$.
Then
$$
\Phi_G^{\ell,i}(\Zym{i})= \ym{2\ell+i-1}{0}.
$$
\end{lm}

There are inclusions $\Zfl{i}\hra\Zym{i}$ given by
$(V,v,\bV,\bv)\mapsto (V,v,\bV,\bv,0)$, where
$V,\bV\in G^{2\ell}$, and $v,\bv\in G^i$.
We use the same notation for $\Phi_G^{\ell,i}$ in Lemma \ref{thm:ZXfl} and Lemma \ref{thm:ZXym},
since $\Phi_G^{\ell,i}$ in Lemma \ref{thm:ZXfl} is just the restriction
of $\Phi_G^{\ell,i}$ in Lemma \ref{thm:ZXym}.

\begin{proof}[Proof of Lemma \ref{thm:ZXym} ]

\noindent
{\bf 1.} Claim: $\Phi_G^{\ell,1}(\Zym{1})\subset \ym{2\ell}{0}$.

Given $(V,c,\bV,\bc,X)\in \Zym{1}$, where $V,\bV\in G^{2\ell}$, $c,\bc\in G$,
and $X\in \fg$, we have
$$
  V,  c\bV c^{-1} \in (G_X)^{2\ell},\quad
  \fm(V)=\exp(X/2)c\bc,\quad
\fm(\bV)=\bc \exp(-X/2)c
$$
Straightforward calculations show that
 $c\bc=\exp(-X/2)\fm(V) \in G_X$,
and
$$\fm(c\fr(\bV)c^{-1})=c\, \fm(\fr(\bV))\,c^{-1}=\exp(X/2) (c\bc)^{-1}.$$
So if $(V,c,\bV,\bc,X)\in \Zym{1}$, then
\begin{eqnarray*}
&& \fm(V)\fm(c\,\fr(\bV)\,c^{-1})
=\exp(X/2)c\bc \exp(X/2)(c\bc)^{-1}\\
&& =\exp(X/2)\exp(\Ad(c\bc)(X)/2)= \exp(X)
\end{eqnarray*}
where we have used $c\bc\in G_X$, i.e.,
$\Ad(c\bc)(X)=X$. In other words,
$$
\Phi_G^{\ell,1}(V,c,\bV,\bc,X)=(V, c\,\fr(\bV) c^{-1},X)
\in \ym{2\ell}{0}.
$$

\noindent
{\bf 2.} Claim: $\ym{2\ell}{0}\subset \Phi_G^{\ell,1}(\Zym{1})$.

Given  $(V_1,V_2,X)\in \ym{2\ell}{0}$, where $V_1,V_2\in G^{2\ell}$ and $X\in \fg$,
we have
$$
V_1,V_2\in (G_X)^{2\ell},\quad\fm(V_1)\fm(V_2)=\exp(X).
$$
Let $\bc=\exp(-X/2)\fm(V_1)$. Then $\bc\in G_X$, and
$$
\bc=\exp(-X/2)\exp(X)\fm(V_2)^{-1}=\exp(X/2)\fm(\fr(V_2))
$$
We have
$$
\fm(V_1)=\exp(X/2)\bc,\quad \fm(\fr(V_2))= \exp(-X/2)\bc=\bc\exp(-X/2),
$$
so
$$
(V_1, e,\fr(V_2),\bc,X)\in \Zym{1},
$$
and
$$
(V_1, V_2,X)=\Phi_G^{\ell,1} (V_1, e, \fr(V_2),\bc,X)\in \Phi_G^{\ell,1}(\Zym{1}).
$$

\noindent
{\bf 3.} Claim: $\Phi^{\ell,2}_G(\Zym{2})\subset \ym{2\ell+1}{0}$.

Given $(V,d,c,\bV,\bd,\bc,X)\in \Zym{2}$, where $V,\bV\in G^{2\ell}$, $d,c,\bd,\bc\in G$,
and $X\in \fg$, we have
$$
V,d^{-1}c\,\fr(\bV)\,c^{-1}d \in (G_X)^{2\ell},\quad
\fm(V)=\exp(X/2)c\bd c^{-1}d,\quad
\fm(\bV)=\bc d\exp(-X/2)\bc^{-1}\bd.
$$
Straightforward computations show that
\begin{eqnarray*}
&& \fm(V)\fm(d^{-1}c\,\fr(\bV)\, c^{-1}d)[d^{-1},c\bc]=\fm(V) d^{-1}c \fm(\bV)^{-1} c^{-1}d[d^{-1},c\bc]\\
&& =\exp(X/2)\exp(\Ad(c\bc)(X)/2) = \exp(X)
\end{eqnarray*}
where we also used $c\bc\in G_X$, i.e.,
$\Ad(c\bc)(X)=X$. In other words,
$$
\Phi_G^{\ell,2}(V,d,c,\bV,\bd,\bc,X)=
(V, d^{-1}c\, \fr(\bV)c^{-1}d, d^{-1},c\bc,X)\in \ym{2\ell+1}{0}.
$$

\noindent
{\bf 4.} Claim: $\ym{2\ell+1}{0} \subset \Phi_G^{\ell,2}(\Zym{2})$.

Given $(V_1,V_2,a,b,X)\in \ym{2\ell+1}{0}$, where
$V_1,V_2\in G^{2\ell}$ and $a,b\in G$, we have
$$
V_1,V_2\in (G_X)^{2\ell},\quad
a,b\in  G_X,\quad
\fm(V_1)\fm(V_2)[a,b]=\exp(X).
$$
Let
$$
d=a^{-1},\quad
\bd=a\exp(-X/2)\fm(V_1),\quad
c=a^{-1},\quad
\bc=a b.
$$
Then
$$
\exp(X/2)c\bd c^{-1} d =\fm(V_1),
$$
and
\begin{eqnarray*}
&& \bc d \exp(-X/2)\bc^{-1} \bd
=  a b a ^{-1} \exp(-X/2) b^{-1} \exp(-X/2)\fm(V_1)\\
&& =  [a, b]\exp(-\Ad(b)(X)/2)\exp(-X/2) \exp(X)(\fm(V_2)[a,b])^{-1}= \fm(\fr(V_2))
\end{eqnarray*}
where we have used that $b\in G_X$ (i.e. $\Ad(b)(X)=X$) in the last equality.
So
$$
(V_1, d,c, \fr(V_2),\bd,\bc,X)\in \Zym{2},
$$
and
$$
(V_1,V_2,a, b,X)=\Phi_G^{\ell,2} (V_1, d,c, \fr(V_2),\bd,\bc,X).
$$
\end{proof}

Let $G^2$ act on $\Zym{1}$, $\Zym{2}$ by
$$
(g_1,g_2)\cdot(V,c,\bV,\bc,X)=
(g_1 V g_1^{-1}, g_1 c g_2^{-1},g_2 \bV g_2^{-1},g_2 \bc g_1^{-1},\Ad(g_1)(X)),
$$
$$
(g_1,g_2)\cdot(V,d,c,\bV,\bd,\bc,X)
= (g_1 V g_1^{-1},g_1 d g_1^{-1}, g_1 c g_2^{-1},g_2 V g_2^{-1},
g_2 \bd g_2^{-1}, g_2 \bc g_1^{-1},\Ad(g_1)(X)),
$$
respectively, where $V,\bV\in G^{2\ell}$, and $g_1,g_2,c,\bc,d,\bd\in G$.
Slight modification of the proof of Lemma \ref{thm:ZGGXG} gives the
following:
\begin{lm}
The surjection $\Phi^{\ell,i}_G: \Zym{i}\to \ym{2\ell+i-1}{0}$
induces homeomorphisms
\begin{equation}
\Zym{i}/ G^2 \cong \ym{2\ell+i-1}{0}/G \cong \Hom(\Gamma_\bR(\Si^{2\ell+i-1}_0),G)/G
\end{equation}
and a homotopy equivalence 
\begin{equation}
\Zym{i}^{h G^2}\sim \ym{2\ell+i-1}{0}^{h G}
\end{equation}
between homotopic orbit spaces.
\end{lm}

\subsection{Involution on representation varieties for Yang-Mills connections}

\begin{lm}\label{thm:VbV}
For $i=1,2$, define
$\tau:G^{2(2\ell+i)}\times\fg \to G^{2(2\ell+i)}\times \fg$ by
\begin{eqnarray*}
\tau(V,c,\bV,\bc,X)&=&(\bV,\bc,V,c,-\Ad(\bc)X)\\
\tau(V,d,c,\bV,\bd,\bc,X)&=& (\bV,\bd,\bc,V,d,c,-\Ad(\bc)X)
\end{eqnarray*}
where $V\in G^{2\ell}$, $c,\bc,d,\bd\in G$, $X\in \fg$. Then
$$
\tau(\Zym{i})=\Zym{i}
$$
and $\tau\circ \tau$ restricts to the identity map on $\Zym{i}$.
\end{lm}

\begin{rem}
Based on \eqref{eqn:inv_1} and \eqref{eqn:inv_2}, the involution $\tau$ defined in
Lemma \ref{thm:VbV} is
the one induced by the $\bZ/2\bZ$ deck transformation on the double cover.
\end{rem}

\begin{proof}[{Proof of Lemma \ref{thm:VbV}}]
We first prove $\tau(\Zym{i})\subset \Zym{i}$.

\noindent
$i=1$:
Given $(V,c,\bV,\bc,X)\in \Zym{1}$, where
$V,\bV\in G^{2\ell}$, $c,\bc\in G$, and $X\in \fg$, we have
$$
V,c\bV c^{-1}\in (G_X)^{2\ell},\quad c\bc \in G_X
$$
so $
\bc^{-1} \bV \bc= (c\bc)^{-1} (c \bV c^{-1} ) (c\bc)
\in (G_X)^{2\ell},
$ or equivalently,
\begin{enumerate}
\item[(i)]$\bV \in (\bc G_X  \bc^{-1})^{2\ell}=
(G_{\Ad(\bc)(X) })^{2\ell}= (G_{-\Ad(\bc)(X) })^{2\ell}$.
\end{enumerate}
If we let
$\bX=-\Ad(\bc)(X)$, then we have $\bV \in (G_{\bX})^{2\ell}$.
We also have
\begin{enumerate}
\item[(ii)]$
\bc V \bc^{-1} \in (\bc G_X \bc^{-1})^{2\ell}
=(G_{\Ad(\bc)(X)})^{2\ell}= (G_{\bX})^{2\ell}
$\end{enumerate}
To summarize, we have
\begin{equation}\label{eqn:RPi}
\bV, \bc V \bc^{-1}\in (G_{\bX})^{2\ell}.
\end{equation}
We also have
\begin{equation}\label{eqn:RPii}
\fm(\bV)= \bc \exp(-X/2) c
=(\bc\exp(-X/2)\bc^{-1})\bc c=\exp(\bX/2) \bc c
\end{equation}
\begin{equation}\label{eqn:RPiii}
\fm(V)=\exp(X/2)c\bc= c\bc \exp(X/2)
=c (\bc \exp(X/2) \bc^{-1})\bc = c \exp(-\bX/2) \bc
\end{equation}
By (\ref{eqn:RPi}), (\ref{eqn:RPii}), (\ref{eqn:RPiii}), we get
$\tau(V,c,\bV,\bc, X)=(\bV,\bc,V,c,\bX) \in \Zym{1}$. This proves
$$
\tau(\Zym{1})\subset \Zym{1}.
$$

\noindent
$i=2$:
Given $(V,d,c,\bV,\bd,\bc,X) \in \Zym{2}$, where
$V,\bV\in G^{2\ell}$, $d,c,\bd,\bc\in G$, and $\fg\in X$, we have
$$
V, d^{-1}c\bV c^{-1}d\in G_X^{2\ell},\quad  c\bc \in G_X,
$$
so $\bc^{-1}\bV \bc=(c\bc)^{-1}d(d^{-1}c\bV c^{-1}d)d^{-1}(c\bc)\in G_X^{2\ell}$, or equivalently,
\begin{enumerate}
\item[(i)]  $ \bV \in (\bc G_X  \bc^{-1})^{2\ell}=
(G_{\Ad(\bc)(X) })^{2\ell}= (G_{\bX})^{2\ell}$,
\end{enumerate}
where $\bX=-\Ad(\bc)(X)$.
We also have
\begin{enumerate}
\item[(ii)] $\bc^{-1}(\bc c)\bc=c\bc\in G_X$, i.e. $\bc c \in G_{\Ad(\bc)(X)}=G_{\bX}$.
\end{enumerate}
On the other hand,
\begin{enumerate}
\item[(iii)] $
\bc^{-1}\bd\bc=(c\bc)^{-1}\exp(-X/2)
\fm(V) d^{-1}(c\bc) \in G_X,\quad \textup{ i.e. }\bd \in G_{\bX}.
$
\end{enumerate}
\begin{enumerate}
\item[(iv)]
$V\in (G_X)^{2\ell}$, so $\bc V \bc^{-1}\in (G_{\Ad(\bc)X})^{2\ell}=(G_{\bX})^{2\ell}$, and thus
$\bd^{-1}\bc V \bc^{-1}\bd \in (G_{\bX})^{2\ell}$.
\end{enumerate}

To summarize, we have
\begin{equation}\label{keqn:i}
V, \bd^{-1}\bc V \bc^{-1}\bd \in (G_X)^{2\ell},\quad
\bd^{-1}, \bc c \in G_{\bX}
\end{equation}
We also have
\begin{equation}\label{keqn:ii}
\fm(\bV)=\bc d\exp(-X/2)\bc^{-1}\bd
=\bc\exp(-X/2) d \bc^{-1}\bd
=\exp(\bX/2)\bc d\bc^{-1}\bd
\end{equation}
\begin{equation}\label{keqn:iii}
\begin{aligned}
\fm(V)& =\exp(X/2)c\bd c^{-1} d
= c\bd\cdot \bd^{-1} (\bc c)^{-1} \bc \exp(X/2)\bc^{-1}
        (\bc c)\cdot \bd \cdot c^{-1}d\\
& =c\bd \exp(\Ad\left(\bd^{-1}\right)\circ\Ad\left((\bc c)^{-1}\right)
\left(-\bX)/2\right) c^{-1} d
=c\bd\exp(- \bX/2)c^{-1}d
\end{aligned}
\end{equation}
where we have used  $\bd^{-1}, (\bc c)^{-1}\in G_{\bX}$ in the last equality.
By (\ref{keqn:i}), (\ref{keqn:ii}), (\ref{keqn:iii}),
$$
\tau(V,d,c,\bV,\bd,\bc,X)=(\bV,\bd,\bc,V,d,c,\bX) \in \Zym{2}.
$$
This proves
$$
\tau(\Zym{2})\subset \Zym{2}.
$$

It remains to show that $\tau\circ \tau:\Zym{i}\to \Zym{i}$
is the identity map. We first consider the case $i=1$:
given $(V,c,\bV,\bc,X)\in\Zym{1}$,
$$
\tau\circ\tau(V,c,\bV,\bc,X)
=\tau(\bV,\bc,V,c,-\Ad(\bc)(X)
=(V,c,\bV,\bc, \Ad(c\bc)(X))
=(V,c,\bV,\bc,X)
$$
where we have used $c\bc\in G_X$. The case
$i=2$ can be proved in the same way.

Thus we have
$$
\tau(\Zym{i})\subset \Zym{i},\quad
\Zym{i}=\tau\circ \tau(\Zym{i})\subset \tau(\Zym{i}).
$$
\end{proof}

\subsection{Representation varieties for nonorientable surfaces}

From the above discussion we have
\begin{equation}\label{eqn:NGZ}
\bigcup_{P\in \Prin_G(\Si^\ell_i)}\cN(P)/\cG(P)
\cong \Zym{i}^\tau/(G\times G)^\tau
\end{equation}
for $i=1,2$, where $\tau:G\times G\to G\times G$ is given by $(g_1,g_2)\mapsto (g_2,g_1)$.
We now relate the right hand side of \eqref{eqn:NGZ} to representations
of the super central extension $\Gamma_\bR(\Si^\ell_i)$ of $\pi_1(\Si^\ell_i)$.

$\Gamma_\bR(\Si^\ell_1)$ is generated by
$$
A_1,B_1,\ldots, A_\ell, B_\ell, C, J_r
$$
where $r\in\bR$, with relations
\begin{enumerate}
\item $J_{r_1}J_{r_2}= J_{r_1+r_2}$
\item $A_i J_r A_i^{-1}=B_i J_r B_i^{-1}=J_r,\quad  i=1,\ldots,\ell,\quad r\in\bR$.
\item $CJ_r C^{-1}=J_{-r}, \quad r\in\bR$.
\item $\prod_{i=1}^\ell [A_i,B_i] = J_1 C^2$.
\end{enumerate}
Let $\rho: \Gamma_\bR(\Si^\ell_1)\to G$ be a group homomorphism.
From the relation (1) we must have
$\rho(J_r)=\exp(rX)$
for some $X\in \fg$. From the relation (2) we must have
$\rho(A_i), \rho(B_i)\in G_X$.
From the relation (3) we have
$\Ad( \rho(C) )(X)=-X$.
Combined with the relation (4), $\Hom(\Gamma_\bR(\Si^\ell_1),G)$ can be identified
with
$$
\ym{\ell}{1}=\{ (V,c,X)\in G^{2\ell+1}\times \fg\mid \\
 V\in (G_X)^{2\ell} , \Ad(c)(X)=-X,  \fm(V)=\exp(X)c^2 \}
$$
There is a homeomorphism
$\ym{\ell}{1} \stackrel{\cong}{\longrightarrow}\Zym{1}^\tau$
given by
$$
(V,c,X)\mapsto (V,c,V,c, 2X),\quad V\in G^{2\ell},c\in G,X\in \fg.
$$
There is an inclusion $\fl{\ell}{1}\hra \ym{\ell}{1}$ given by
$(V,c)\mapsto (V,c,0)$.

$\Gamma_\bR(\Si^{\ell}_2)$ is generated by
$$
A_1,B_1,\ldots, A_\ell, B_\ell, D, C, J_r
$$
where $r\in\bR$, with relations
\begin{enumerate}
\item $J_{r_1}J_{r_2}= J_{r_1+r_2}$
\item $A_i J_r A_i^{-1}=B_i J_r B_i^{-1}=D J_r D^{-1}=J_r,\quad  i=1,\ldots,\ell,\quad r\in\bR$.
\item $CJ_r C^{-1}=J_{-r}, \quad r\in\bR$.
\item $\prod_{i=1}^\ell [A_i,B_i] = J_1 CDC^{-1}D$.
\end{enumerate}

Let $\rho: \Gamma_\bR(\Si^\ell_2)\to G$ be a group homomorphism.
From the relation (1) we must have $\rho(J_r)=\exp(rX)$
for some $X\in \fg$. From
the relation (2) we must have
$\rho(A_i), \rho(B_i), \rho(D)\in G_X$.
From the relation (3) we have
$\Ad( \rho(C) )(X)=-X$.
Combined with the relation (4), $\Hom(\Gamma_\bR(\Si^\ell_2),G)$ can be identified
with
\begin{eqnarray*}
\ym{\ell}{2}&=& \{ (V, d, c, X)\in G^{2\ell+2}\times \fg\mid
V\in (G_X)^{2\ell}, d\in G_X,\\
&& \quad  \Ad(c)(X)=-X,\ \fm(V)=\exp(X)cdc^{-1}d\}
\end{eqnarray*}
There is a homeomorphism
$\ym{\ell}{2} \stackrel{\cong}{\longrightarrow}\Zym{2}^\tau$
given by
$$
(V,d,c,X)\mapsto (V,d,c,V,d,c, 2X),\quad V\in G^{2\ell},d,c\in G,X\in \fg.
$$
There is an inclusion $\fl{\ell}{2} \hra  \ym{\ell}{2}$ given by
$(V,d,c)\mapsto (V,d,c,0)$.

We obtain the following analogue of Theorem \ref{thm:YMzero} for nonorientable surfaces.
\begin{thm}\label{thm:YMonetwo}
There is a bijective correspondence between conjugacy classes of homomorphisms
$\Gamma_\bR(\Si)\to G$ and gauge equivalence classes of Yang-Mills $G$-connections over $\Si$.
In other words, for $i=1,2$, we have
\begin{eqnarray*}
\bigcup_{P\in \Prin_G(\Si^\ell_i)} \cN(P)/\cG_0(P)
&\cong& \Hom(\Gamma_\bR(\Si^\ell_i),G)\cong \ym{\ell}{i},\\
\bigcup_{P\in \Prin_G(\Si^\ell_i)} \cN(P)/\cG(P)
&\cong&\Hom(\Gamma_\bR(\Si^\ell_i),G)/G \cong \ym{\ell}{i}/G,
\end{eqnarray*}
where $g\in G$ acts on $G^{2\ell+1}\times \fg$ by
$$
g\cdot (V,c,X)= (gVg^{-1},g c g^{-1}, \Ad(g)(X) ),
$$
and on $G^{2\ell+2}\times \fg$ by
$$
g\cdot (V,d,c,X)=(gVg^{-1}, gdg^{-1},gcg^{-1},\Ad(g)(X)).
$$
\end{thm}

\subsection{Extended moduli spaces}
The representation variety $\ym{\ell}{0}$ is a subset of Lisa Jeffrey's
extended moduli space \cite{J}.
In this subsection, we define extended moduli
spaces for nonorientable surfaces.

Let $\Si^{\ell,r}_0$ be the compact, connected, orientable surface with $\ell$ handles and
$r$ boundary components $S_1,\ldots,S_r$ with coordinates $(s_1,\ldots,s_r)\in \bR/\bZ$.
Let $\Si^{\ell,r}_1$ be the connected sum of $\Si^\ell_0$ and $\RP^2$, and let
$\Si^{\ell,r}_2$ be the connected sum of $\Si^\ell_0$ and a Klein bottle.
The following discussion is a straightforward
generalization of the case $i=0$ in \cite{J}.

Suppose that $r>0$ and $i=0,1,2$. Then any principal $G$-bundle
$P$ over $\Si^{\ell,r}_i$ is topologically trivial. Let
$\cA(P)$ be the space of smooth connections on $P$. Then
$$
\cA(P)\cong \Omega^1_{\Si^{\ell,r}_i}(\fg).
$$
Define
\begin{eqnarray*}
\cA^\fg_G(\Si^{\ell,r}_i)&=& \{ A\in \cA(P)\mid F(A)=0, A|_{U_j}=X_j ds_j
\textup{ on some open neighborhood }\\
&&  U_j \textup{ of } S_j \textup{ for
some }X_j\in \fg, j=1,\ldots,r\},
\end{eqnarray*}
and define the compactly supported gauge group
$$
\cG^c(\Si^{\ell,r}_i)= \{ s:\Si^{\ell,r}_i \to G \mid s(x) = e
\textup{ for }x\in U, \textup{ where } U \textup{ is an open neighborhood of }
\partial \Si^{\ell,r}_i  \}.
$$
We define a moduli space
$$
\cM^\fg_G(\Si^{\ell,r}_i)= \cA^\fg_F(\Si^{\ell,r}_0)/\cG^c(\Si^{\ell,r}_i),
$$
and introduce representation varieties
\begin{eqnarray*}
\cN^\fg_G(\Si^{\ell,r}_0)&=&\{(V, k_2,\ldots,k_r, X_1,\ldots,X_r
\in G^{2\ell}\times G^{r-1}\times \fg^r \mid \\
 && \fm(V) =\exp(X_1)\exp(\Ad(k_2)X_2)\cdots \exp(\Ad(k_r)X_r) \},\\
\cN^\fg_G(\Si^{\ell,r}_1)&=&\{(V,c, k_2,\ldots,k_r, X_1,\ldots,X_r
\in G^{2\ell}\times G^r\times \fg^r \mid \\
 && \fm(V) =\exp(X_1)\exp(\Ad(k_2)X_2)\cdots \exp(\Ad(k_r)X_r) c^2\},\\
\cN^\fg_G(\Si^{\ell,r}_2)&=&\{(V,d,c, k_2,\ldots,k_r, X_1,\ldots,X_r
\in G^{2\ell}\times G^{r+1}\times \fg^r \mid \\
 && \fm(V)=\exp(X_1)\exp(\Ad(k_2)X_2))\cdots \exp(\Ad(k_r)X_r)cdc^{-1}d \}.
\end{eqnarray*}
where $\fm(\ab)=\pab$ as before.
In particular,
\begin{eqnarray*}
\cN^\fg_G(\Si^{\ell,1}_0)&=&\{(V,  X )\in G^{2\ell}\times \fg \mid \fm(V) =\exp(X) \},\\
\cN^\fg_G(\Si^{\ell,1}_1)&=&\{(V,c, X) \in G^{2\ell+1}\times \fg \mid \fm(V) =\exp(X) c^2\},\\
\cN^\fg_G(\Si^{\ell,1}_2)&=&\{(V,d,c, X) \in G^{2\ell+2}\times \fg \mid \fm(V) =\exp(X)cdc^{-1}d \},
\end{eqnarray*}
so
$$
\ym{\ell}{i}\subset \cN^\fg_G(\Si^{\ell,1}_i).
$$

The following statement follows from the proof of \cite[Proposition 5.3]{J}:
\begin{pro}Let $\ell\geq 0, r>0$ be integers, and let $i=0,1,2$.
Then there is a homeomorphism
$$
\cM^\fg_G(\Si^{\ell,r}_i)\cong \cN_G^\fg(\Si^{\ell,r}_i).
$$
\end{pro}

\section{Equivariant Morse Stratification of Space of Connections}
\label{sec:general}

In this section, we discuss the $\cG$-equivariant
Morse stratification and reduction
for general compact connected Lie groups.
As a byproduct, we reproduce and extend the results
on connected components of the moduli space of flat connections over closed
nonorientable surfaces \cite{HL2,HL3}.

\subsection{Morse stratification with involution}\label{sec:generalMB}

Let $(M,g)$ be a Riemannian manifold. Let $f:M\to \bR$ be
a smooth function, and let $\phi_t$ be the gradient flow of $f$.
Suppose that the gradient flow is defined for any time $t\in\mathbb{R}$
and the limits
$$
\lim_{t\to \infty}\phi_t(x),\quad
\lim_{t\to -\infty}\phi_t(x)
$$
exist for any $x\in M$. Let $N$ be the set of critical points of $f$, and let
$$
N=\bigcup_{\mu\in\Lambda} N_\mu
$$
be the union of connected components. Suppose that each $N_\mu$ is a closed nondegenerate
critical submanifold of $M$. Given a critical subset $N_\mu$, define its stable manifold $S_\mu$
and unstable manifold $U_\mu$ by
$$
S_\mu =  \{ x\in M\mid \lim_{t\to-\infty}\phi_t(x)\in N_\mu\},\quad
U_\mu = \{x\in M\mid \lim_{t\to +\infty}\phi_t(x)\in N_\mu\}.
$$
Then
$$
M=\bigcup_{\mu\in\Lambda} S_\mu
$$
is a disjoint union of Morse strata. We assume that each $S_\mu$
is a submanifold of $M$.

Suppose that $f$ is invariant under some isometric involution
$\tau: M\to M$. By Lemma \ref{thm:grad}, $\tau$ induces an involution
$\tau_0:\Lambda\to \Lambda$ such that
$$
\tau(N_\mu)=N_{\tau_0(\mu)},\quad
\tau(S_\mu)=S_{\tau_0(\mu)},\quad
\tau(U_\mu)=U_{\tau_0(\mu)}.
$$

\begin{pro}\label{thm:J}
Let $(M,\omega,J)$ be an almost K\"{a}hler manifold
with  an anti-symplectic, anti-holomorphic involution $\tau:M\to M$.
Suppose that $f:M\to \bR$ is  a $\tau$-invariant smooth function.
Suppose that $N_\mu$ is a closed subset of M and a connected component
of the set of critical points $N$ of $f$. Suppose that the set
$$
S_\mu =\{ x\in M\mid \lim_{t\to -\infty}\phi_t(x)\in N_\mu \}
$$
is an almost complex submanifold of $M$. If
$$
S_\mu^\tau=\{ x\in M^\tau \mid \lim_{t\to -\infty}\phi_t(x)\in N^\tau_\mu \}
$$
is nonempty, then $S_\mu^\tau$ is the stable manifold of $\cN^\tau_\mu$ with
respect to $f^\tau$, and the real codimension
of $S_\mu^\tau$ in $M^\tau$ is equal to the complex codimension
of $S_\mu$ in $M$.
\end{pro}

Note that $M^\tau$ and $S^\tau_\mu$ are not necessarily connected.

\subsection{Morse stratification and Morse inequalities}
\label{sec:Mstrata}

Let $\Si$ be a closed, compact, connected surface.
Then $\Si$ is diffeomorphic to $\Si^\ell_i$ for some integer $\ell\geq 0$
and $i\in\{0,1,2\}$. Recall that $\chi(\Si^\ell_i)= 2-2\ell -i$.

Let $G$ be a compact, connected Lie group. By Theorem \ref{thm:YMzero}
and Theorem \ref{thm:YMonetwo},
$$
\bigcup_{P\in \Prin_G(\Si)} \cN(P)/\cG_0(P) \cong \Hom(\Gamma_\bR(\Si),G),\quad
\bigcup_{P\in \Prin_G(\Si)} \cN(P)/\cG(P) \cong \Hom(\Gamma_\bR(\Si),G)/G.
$$
Let
$$
\Hom(\Gamma_\bR(\Si),G)_P \subset \Hom(\Gamma_\bR(\Si),G)
$$
be the subset corresponding to $P\in \Prin_G(\Si)$, so that
$$
\Hom(\Gamma_\bR(\Si),G)=\bigcup_{P\in\Prin_G(\Si)} \Hom(\Gamma_\bR(\Si),G)_P,
$$
$$
\cN(P)/\cG_0(P) \cong \Hom(\Gamma_\bR(\Si),G)_P,\quad
\cN(P)/\cG(P) \cong \Hom(\Gamma_\bR(\Si),G)_P/G.
$$

For a fixed topological principal $G$-bundle $P$ over $\Si$, let
$$
\{ \cN_{\tilde{\mu}} (P) \mid \tilde{\mu}\in \pi_0(\cN(P))\}
$$
be the set of connected components of $\cN(P)$.  
Let $\cG(P)'$ be the connected component of
the identity of $\cG(P)$. Then $\cG(P)'$ is 
a normal subgroup of $\cG(P)$, and the 
discrete set $\pi_0(\cG(P))$ can be identified
with the group $\cG(P)/\cG(P)'$. The action of $\cG(P)$
on $\cN(P)$ is continuous, and induces an action
of $\pi_0(\cG(P))\cong \cG(P)/\cG(P)'$ on 
$\pi_0(\cN(P))$. Define
$$
\Lambda = \pi_0(\cN(P))/\pi_0(\cG(P)).
$$

\begin{rem}\label{pi-zero}
When $\Si$ is orientable, the action of $\pi_0(\cG(P))$
on $\pi_0(\cN(P))$ is trivial by the results in \cite{ym}, so
$\Lambda = \pi_0(\cN(P))$. 
\end{rem}

Let $\pi:\pi_0(\cN(P))\to \Lambda=\pi_0(\cN(P))/\pi_0(\cG(P))$ be the projection.
Given $\mu\in \Lambda$, define
$$
\cN_\mu(P)=\bigcup_{\tilde{\mu}\in \pi^{-1}(\mu)}\cN_{\tilde{\mu}}(P).
$$
Then $\cG(P)$ acts on $\cN_\mu(P)$, and
$\cN_\mu(P)/\cG(P)$ is connected. Note that
$\cG(P)/\cG_0(P)=G$ is connected, so
$$
V_\mu(P)\define \cN_\mu(P)/\cG_0(P)
$$
is connected.
$$
\Hom(\Gamma_\bR(\Si),G)_P =\bigcup_{\mu\in \Lambda} V_\mu(P)
$$
is a disjoint union of connected components.
Each $\cN_\mu(P)$ is a closed subset of $\cN(P)$ thus of $\cA(P)$. Define
$$
\cA_\mu(P)=\{ A\in \cA(P)\mid \lim_{t\to -\infty} \phi_t(A)\in \cN_\mu(P) \}
$$
where $\phi_t$ is the gradient flow of $L_P$.
The limit exists by results in \cite{da} and \cite{Rad}. 
Notice that $L_P$ is constant on each
 $\cN_\mu(P)$, and $L_P$ achieves its minimum on $\cN_\mu(P)$ within $\cA_\mu(P)$. 
Each $\cA_\mu(P)$ is
a submanifold of $\cA(P)$, and the map $\cA_\mu(P) \rightarrow \cN_\mu(P)$ given
by $A\mapsto\displaystyle{\lim_{t\to -\infty} \phi_t(A)}$ is a $\cG(P)$-equivariant
deformation retraction. We have
\begin{equation}\label{eqn:strataG}
\cA(P)=\bigcup_{\mu\in\Lambda}\cA_\mu(P)
\end{equation}
is a smooth stratification. This stratification
is $\cG(P)$-equivariant in the sense
that $\cG(P)$ acts on each stratum.
We call \eqref{eqn:strataG} the  $\cG(P)$-equivariant Morse stratification
of $\cA(P)$ with respect to $L_P$. 

\begin{rem}
Given $\tilde{\mu}\in \pi_0(\cN(P))$, define
$\cA_{\tilde{\mu}}(P)=\{ A\in \cA(P)\mid
\displaystyle{\lim_{t\to -\infty} \phi_t(A)}\in \cN_{\tilde{\mu}}(P) \}$.
Then $\cA_{\tilde{\mu}}$ is a connected submanifold of $\cA(P)$, and is a Morse stratum
of $L_P$. When $\Si$ is orientable, the Morse stratification
coincides with the $\cG(P)$-equivariant Morse stratification
by Remark \ref{pi-zero}; when $\Si$ is nonorientable, a priori the Morse stratification
can be finer than the $\cG(P)$-equivariant Morse stratification.
 \end{rem}

We now assume that $\chi(\Si)<0$.
Let $\cN_{ss}$ be the set where the Yang-Mills functional $L_P$ achieves
absolute minimum.   Then $\cN_{ss}$ is connected
(by results in \cite{ym} when $\Si$ is orientable, and
by Theorem \ref{thm:Nconnected} 
when $\Si$ is nonorientable). 
$\cN_{ss}$ is the unique connected component of
the critical set $\cN(P)$ with zero Morse index, and its
stable manifold $\cA_{ss}$ is the unique codimension
zero Morse stratum
(which is also a $\cG(P)$-equivariant Morse stratum). 
Define
$$
V_{ss}(P)= \cN_{ss}/\cG_0(P).
$$
Then $V_{ss}(P)$ is connected. 

When the obstruction class $o(P)\in H^2(\Si,\pi_1(G))$
is a torsion element (which is always true when $\Si$ is nonorientable),
$\cN_{ss}(P)=\cN_0(P)$ is the space of flat connections on $P$,
and
$$
V_{ss}(P)\cong \Hom(\pi_1(\Si),G)_P
$$
where $\Hom(\pi_1(\Si),G)_P$
is the connected component of $\Hom(\pi_1(\Si),G)$
associated to the topological $G$-bundle $P$ (see
\cite{ym, HL3}).

In general, we are interested in the cohomology of the moduli space
$\cM(P)$ of gauge equivalence classes of minimal Yang-Mills connections
on $P$. More explicitly,
$$
\cM(P)\define \cN_{ss}(P)/\cG(P) \cong V_{ss}(P)/G.
$$
When $\cM(P)$ is smooth and $\cG(P)$ acts on $\cN_{ss}(P)$ freely, we have
$$
H^*(\cM(P);\bQ)\cong H^*_{\cG(P)}(\cN_{ss}(P);\bQ)\cong H^*_G(V_{ss}(P);\bQ).
$$

The deformation retraction $r:\cA_\mu(P)\to \cN_\mu(P)$ given by the
gradient flow of $L_P$ is $\cG(P)$-equivariant, thus the following equivariant pairs
are equivalent for the purpose of (singular) equivariant cohomology:
$$
(\cA_\mu(P),\cG(P))\sim (\cN_\mu(P),\cG(P))\sim (\cN_\mu(P)/\cG_0(P), \cG(P)/\cG_0(P))
\sim(V_\mu(P),G).
$$
In other words, we have the following homotopy equivalences of homotopic orbit spaces:
$$
\cA_\mu(P)^{h\cG(P)} \sim \cN_\mu(P)^{h \cG(P)} \sim V_\mu(P)^{h G}.
$$
As a consequence, we have the following isomorphisms of (singular) equivariant cohomology:
\begin{equation}\label{eqn:ANV}
H^*_{\cG(P)}(\cA_\mu(P);\bQ) \cong H^*_{\cG(P)}(\cN_\mu(P);\bQ) \cong H^*_G(V_\mu(P);\bQ).
\end{equation}

Let $K$ be a field and let
$$
P^{\cG(P)}_t(\cA_\mu(P);K)=\sum_{i\geq 0} t^i \mbox{dim}H^i_{\cG(P)}(\cA_\mu(P);K)
$$
be the equivariant Poincar\'{e} series.
Let
\begin{equation}\label{eqn:MtcG}
M^{\cG(P)}_t(L_P;K)=\sum_{\mu\in \Lambda}t^{\lambda_\mu} P^{\cG(P)}_t(\cA_\mu(P);K)=
P^{\cG(P)}_t(\cA_{ss}(P);K)+\sum_{\mu\in \Lambda'}t^{\lambda_\mu}P^{\cG(P)}_t(\cA_\mu(P);K)
\end{equation}
be the $\cG(P)$-equivariant Morse series of $L_P:\cA(P)\to \bR$,
where $\lambda_\mu$ is the real codimension of the stratum $\cA_\mu(P)$ in $\cA(P)$
and $\Lambda'=\{\mu\in \Lambda\mid \lambda_\mu>0 \}$.

The Morse stratification is smooth and $\cG(P)$-equivariant, so we have equivariant Morse inequalities
\cite{da, K1}: there
exists a power series $R_K(t)$ with nonnegative coefficients such that
\begin{equation}\label{eqn:inequality}
M^{\cG(P)}_t(L_P;K) -P^{\cG(P)}_t(\cA(P);K)=(1+t)R_K(t).
\end{equation}
A priori (\ref{eqn:inequality}) holds for $K=\bZ/2\bZ$. If the normal bundle
$\nu(\cA_\mu)$ of $\cA_\mu(P)$ in $\cA(P)$ is {\em orientable} then
\eqref{eqn:inequality} holds for any $K$. When $\Si$ is orientable,
$\nu(\cA_\mu)$ is a complex vector bundle and has a canonical orientation.

\eqref{eqn:ANV} and \eqref{eqn:MtcG} imply
\begin{equation}\label{eqn:MtG}
M^{\cG(P)}_t(L_P;K)=\sum_{\mu\in \Lambda}t^{\lambda_\mu} P^G_t(V_\mu(P);K)=
P^G_t(V_{ss}(P);K)+\sum_{\mu\in \Lambda'}t^{\lambda_\mu} P^G_t(V_\mu(P);K).
\end{equation}
where now we consider the equivariant cohomology of $V_\mu(P)$, the representation variety, which is
finite dimensional but singular, as opposed to $\cA_\mu(P)$, which is smooth but infinite dimensional.

\subsection{Equivariant Poincar\'{e} series}
When $\Si$ is orientable, Atiyah and Bott \cite{ym} provide an algorithm to
compute the equivariant Poincar\'{e} series
$$
 P^{\cG(P)}_t(\cA_{ss}(P);\bQ)=P^{\cG(P)}_t(\cN_{ss}(P);\bQ).
$$
We now outline this algorithm.

Let $G^\bC$ be the complexification of $G$. Then
$G^\bC$ is a connected reductive algebraic group over $\bC$.
(For example, if $G=U(n)$ then $G^\bC=GL(n,\bC)$.)
There is a bijection
$$
\Prin_G(\Si)\to \Prin_{G^\bC}(\Si)
$$
given by $P\mapsto P\times_G G^\bC$.

We can consider $\cC(\xi)$, the space of $(0,1)$-connections on $\xi=P\times_G G^\bC$.
Using Harder-Narasimhan filtration, Atiyah-Bott gave a stratification
for $\cC(\xi)$ by $$\cC(\xi)=\bigcup_{\mu\in \Xi}\cC_\mu(\xi)$$ where $\mu$ denotes
 Harder-Narasimha type $\mu$ and each $\cC_\mu$ is connected (\cite[Chapter 7, 10]{ym}) .

Atiyah-Bott proved that the Harder-Narasimhan 
stratification is $\cG^\bC$-equivariantly perfect over $\bQ$, i.e.
\begin{equation}\label{eqn:cGPt}
P^{\cG^\bC}_t(\cC(\xi);\bQ)=P^{\cG^\bC}_t(\cC_{ss}(\xi);\bQ)+
\sum_{\mu\in \Xi'}t^{\lambda_\mu}P^{\cG^\bC}_t(\cC_\mu(\xi);\bQ)
\end{equation}
where $\cC_{ss}$ is the semistable stratum,
$\Xi'=\{\mu\in \Xi\mid \cC_\mu(\xi)\neq \cC_{ss}(\xi)\}$.
Now $\cC(\xi)$ is contractible, so
\begin{equation}\label{eqn:BcG}
H^*_{\cG^\bC}(\cC(\xi);\bZ)= H^*_{\cG^\bC}(\mathrm{pt};\bZ) = H^*(B\cG^\bC;\bZ)=H^*(B\cG(P);\bZ).
\end{equation}
Thus\begin{equation}
P_t(B\cG(P);\bQ)= P^{\cG^\bC}_t(\cC_{ss}(\xi);\bQ)+
\sum_{\mu\in \Xi'}t^{\lambda_\mu} P^{\cG^\bC}_t(\cC_\mu(\xi);\bQ).
\end{equation}

On the other hand, there is a natural isomorphism $i: \cC(\xi)\rightarrow \cA(P)$, and
 it was proven in \cite{da} (conjectured by \cite{ym}) that the Harder-Narasimhan stratification coincides with the Morse stratification defined by $L_P$ as in Section \ref{sec:Mstrata}, i.e. 
 $\Xi\cong\pi_0(\cN(P))\cong \Lambda$.
So the codimension of $\cC_\mu$ in $\cC$
equals to the codimension of $\cA_\mu=i(\cC_\mu)$ in $\cA$. 
In particular, $\cA_{ss}=i(\cC_{ss})$.
The moduli space $\cM(P)$ of minimal (central) Yang-Mills connections on $P$ 
can be identified with the moduli space of $S$-equivalence
classes of semi-stable holomorphic structures on $\xi$ \cite{Ra,ym}.

The equivariant perfectness of Harder-Narasimhan stratification $\bigcup \cC_\mu$ now implies that
\[P^{\cG}_t(\cA(P);\bQ)=P^{\cG^\bC}_t(\cC(\xi);\bQ)=\sum_{\mu\in\Xi}t^{\lambda_\mu} P^{\cG^\bC}_t(\cC_\mu(\xi);\bQ)
=\sum_{\mu\in\Lambda}t^{\lambda_\mu} P^{\cG}_t(\cA_\mu(P);\bQ)=M^{\cG(P)}_t(L_P;\bQ)\]
i.e.  $R_\bQ(t)=0$ and Morse stratification is $\cG(P)$-equivariantly perfect:
\begin{equation}\label{eqn:GPt}
P_t(B\cG(P);\bQ)= P^G_t(V_{ss}(P);\bQ)+
\sum_{\mu\in \Lambda'}t^{\lambda_\mu} P^G_t(V_\mu(P);\bQ).
\end{equation}

The equivariant pair $(V_\mu(P),G)$ can be reduced further (cf: \cite[Section 10]{ym}):
$$
(V_\mu(P),G)\sim (V_{ss}(P_\mu), G_\mu)
$$
where $G_\mu$ is a compact Lie subgroup of $G$ with $\dim_\bR G_\mu < \dim_\bR G$,
and $P_\mu$ is  a principal $G_\mu$-bundle. So we have $P^G_t(V_\mu(P);\bQ)=P^{G_\mu}_t(V_{ss}(P_\mu);\bQ)$ and
\begin{equation}\label{eqn:GPt_ss}
P_t(B\cG(P);\bQ)= P^G_t(V_{ss}(P);\bQ)+
\sum_{\mu\in \Lambda'}t^{\lambda_\mu}P^{G_\mu}_t(V_{ss}(P_\mu);\bQ).
\end{equation}
The left hand side of \eqref{eqn:GPt_ss} can be computed;
a formula for the case $G=U(n)$ is given by \cite[Theorem 2.15]{ym}.
So $P^G_t(V_{ss}(P);\bQ)$ can be computed recursively. The case
$G=U(n)$ is particularly nice because
\begin{enumerate}
\item[U1.] $G_\mu$ is of the form $U(n_1)\times\cdots\times  U(n_r)$ so the induction reduction
does not use other groups.
\item[U2.] $H^*(B\cG(P);\bZ)$ is torsion free, so $\rank H^i(B\cG(P);\bZ)=\dim_{\bQ} H^i(B\cG(P);\bQ)$.
\end{enumerate}
Neither U1 nor U2 is true for a general compact connected Lie group
$G$.

Finally, we point out the difficulties in generalizing the above
approach to nonorientable surfaces.
\begin{enumerate}
\item[N1.] It is not clear if $\nu(\cA_\mu)$ is orientable in general,
so a priori Morse inequalities hold only for $\bZ/2\bZ$:
\begin{equation}\label{eqn:BGRmodtwo}
P_t(B\cG(P);\bZ/2\bZ)= P^G_t(V_{ss}(P);\bZ/2\bZ)+
\sum_{\mu\in \Lambda'}t^{\lambda_\mu} P^G_t(V_\mu(P);\bZ/2\bZ)-(1+t)R_{\bZ/2\bZ}(t).
\end{equation}
\item[N2.] The left hand side of \eqref{eqn:BGRmodtwo} is difficult to compute
when $H^*(B\cG(P);\bZ)$ has 2-torsion elements. When $\Si$ is nonorientable,
$H^*(B\cG(P);\bZ)$ has 2-torsion elements even when $G=U(n)$. (We thank
Paul Selick for pointing this out to us.)
\item[N3.] Suppose that for a particular $P$ we can prove that
$\nu(\cA_\mu)$ is orientable for all $\mu\in \Lambda$.
We can consider rational cohomology
\begin{equation}\label{eqn:BGRQ}
P_t(B\cG(P);\bQ)= P^G_t(V_{ss}(P);\bQ)+
\sum_{\mu\in \Lambda'}t^{\lambda_\mu} P^G_t(V_\mu(P);\bQ)-(1+t)R_{\bQ}(t).
\end{equation}
Although $P_t(B\cG(P);\bQ)$ is easier to handle than $P_t(B\cG(P);\bZ/2\bZ)$,
it is tricky to compute $R_{\bQ}(t)$, which is not necessarily zero.
\item[N4.] When the base $\Si$ of the principal $G$-bundle $P$ is
nonorientable, we still have reduction $(V_\mu(P),G)\sim (V_\mu'(P'),G_\mu)$
where $G_\mu$ is a compact Lie subgroup of $G$ with $\dim_\bR G_\mu< \dim_\bR G$,
but $V_\mu'(P')$ is not of the
form $V_{ss}(P_\mu)$ where $P_\mu$ a principal $G_\mu$-bundle over $\Si$.
\end{enumerate}

We will describe equivariant Morse stratification
and the reduction N4 for $G=U(n)$ in Section \ref{sec:Unonorientable} of this paper,
and for other classical groups in \cite{HL5}. We will address N1--N3 in future
works.

\subsection{Connected components of moduli spaces of flat connections}\label{sec:components}
Let $G$ be any compact connected Lie group.
Let $\Si$ be a closed, compact, connected, nonorientable surface, and let
$\pi:\tSi\to \Si$ be its orientable double cover. Let $P\to \Si$ be
a principal $G$-bundle. 
By Proposition \ref{thm:top_trivial}, the pull back principal $G$-bundle
$\tP=\pi^* P\to \tSi$ is topologically trivial, and there is an
involution $\tilde{\tau}:\tP\to \tP$ which covers the deck transformation
$\tau:\tSi\to \tSi$ such that $\tP/\tilde{\tau}=P$. 

The involution $\tilde{\tau}$ induces an involution $\tilde{\tau}^*$ on the space $\cA(\tP)$
of connections on $\tP$, and $\cA(P)$, the space of connections
on $P$, can be identified with the fixed locus $\cA(\tP)^{\tilde{\tau}^*}$.
Each $\cG(P)$-equivariant Morse stratum
$\cA_\mu(P)$ of $L_P$ is a union of connected components of
$\cA_{\mu}(\tP)^{\tilde{\tau}^*}=\cA_{\mu}(\tP)\cap \cA(\tP)^{\tilde{\tau}^*}$,
where $\cA_{\mu}(\tP)$ is a Morse stratum in $\cA(\tP)$. 
The real codimension of $\cA_\mu(P)$ in $\cA(P)$ is equal to the complex
codimension $d_\mu$ of $\cA_{\mu}(\tP)$ in $\cA(\tP)$ (see Proposition \ref{thm:J}).

Let $\cN_0(\tP)$ and $\cN_0(P)$ be the space of flat connections
on $\tP$ and on $P$, respectively. Then
$\cN_0(P)=\cN_0(\tP)^{\tilde{\tau}^*}$. By discussion
in Section \ref{sec:tau_P}, $\cN_0(P)$ is nonempty.
Let $\cA_{ss}(P)= \cA_{ss}(\tP)\cap \cA(\tP)^{\tilde{\tau}^*}$
be the stable manifold of $\cN_0(P)=\cN_{ss}(P)$, so that it is 
the union of all codimension zero Morse strata.
By results in \cite{ym}, $\cA_{ss}(\tP)$ is connected  
when $\chi(\tSi)<0$.
We will show that $\cA_{ss}(P)$ is connected 
when $\chi(\Si)<0$.

\begin{pro}\label{thm:path}
Given two points $A_0,A_1\in \cA_{ss}(P)$, there exists
a smooth map $\gamma:[0,1]\to \cA(P)$ such that
$\gamma(0)=A_0$, $\gamma(1)=A_1$, and $\gamma$ is transversal
to $\cA_\mu(P)$ if $d_\mu >0$.  In particular,
$\gamma^{-1}(\cA_\mu)$  is empty if $d_\mu>1$.
\end{pro}
\begin{proof}
Let $\Omega=\Omega^1(\Si,\ad P)$ be the vector space associated to the affine space $\cA(P)$.
Given $A_0, A_1 \in\cA_{ss}(P)$, define
$$
\Phi: [0,1]\times \Omega \to \cA(P),\quad
\Phi(t,a)= (1-t)A_0+t A_1 +\sin (\pi t)a.
$$
Note that
$$
\Phi(\frac{1}{2},a)=\frac{1}{2}(A_0+A_1)+ a
$$
so $\Phi$ is surjective.
$d\Phi_{(t,a)}: \bR \times \Omega \to \Omega$ is given by
\begin{equation}\label{eqn:dPhi}
(u,b)\mapsto  \left((A_1-A_0)+\pi \cos(\pi t)a\right) u + \sin(\pi t)b
\end{equation}
where $u\in\bR$, $b\in \Omega$. Given $a\in \Omega$ let
$\gamma_a(t)=\Phi(t,a)$. Then
$$
d\Phi_{(t,a)}(u,b)=d(\gamma_a)_t(u) + \sin(\pi t)b
$$

We claim that $\Phi$ is transversal to $\cA_\mu(P)$ for
any $\mu\in \Lambda'$. Fix $\lambda\in \Lambda'$, we need to
show that
$$
\mathrm{Im} \left(d\Phi_{(t,a)}\right) + T_{\Phi(t,a)}\cA_\mu(P)=T_{\Phi(t,a)}\cA(P)=\Omega
$$
for any $(t,a)\in \Phi^{-1}(\cA_\mu(P))$. Note that
$\Phi(0,a)=A_0 \in \cA_{ss}$ and $\Phi(1,a)=A_1 \in \cA_{ss}(P)$ for
any $a\in \Omega$, so if $(t,a)\in\Phi^{-1}(\cA_\mu(P))$ we must
have $0<t<1$. By \eqref{eqn:dPhi}, $\mathrm{Im} \left(d\Phi_{(t,a)}\right)=\Omega$
if $0<t<1$. So $\Phi$ is transversal to $\cA_\mu(P)$ for any $\mu\in \Lambda'$.
We conclude that $\Phi^{-1}(\cA_\mu(P))$ is a submanifold of $[0,1]\times \Omega$ of codimension
$d_\mu$; it is nonempty because $\Phi$ is surjective.

For any $\mu\in \Lambda'$, we define
$\pi_\mu :\Phi^{-1}(\cA_\mu)\to \Omega$ by $(t,a)\mapsto a$.
By Lemma \ref{thm:XYZ},  $a$ is a regular value of $\pi_\mu$ if and only of $\gamma_a:[0,1]\to \cA(P)$
is transversal to $\cA_\mu(P)$. Let $\Omega_\mu$ be the set of regular
values of $\pi_\mu$. By Sard-Smale theorem, $\Omega_\mu$ is residual in $\Omega$. So
$$
\Omega'=\bigcap_{\mu\in \Lambda'}\Omega_\mu
$$
is residual in $\Omega$. By Baire category theorem, $\Omega'$ is nonempty.
For any $a\in \Omega'$, $\gamma_a:[0,1]\to \cA(P)$ has
the desired properties.
\end{proof}

\begin{lm}\label{thm:XYZ}
Let $X,Y,Z$ be linear spaces, and let $W$ be a linear subspace of $Z$.
Let $L: X\times  Y \to Z$ be a linear map such that
$\mathrm{Im}(L)+ W =Z$. Let $\pi: L^{-1}(W)\to Y$ be defined
by $(x,y)\mapsto y$. Then $\pi$ is surjective if and only if
$L(X\times \{0\})+ W =Z$.
\end{lm}
\begin{proof}

\begin{enumerate}
\item Claim: $\pi$ is surjective $\Rightarrow$ $L(X\times \{0\})+W =Z$.

Given any $z\in Z$, we have $z = L(x,y) + w$ for some $(x,y)\in X\times Y$
and $w\in W$. Since $\pi$ is surjective, there is $(x',y')\in X\times Y$
such that $L(x',y') = w'\in W$ and $\pi(x',y')=y$. We have
$$
L(x,y)+ w =z,\quad L(x',y')-w'=0,\quad y=y',
$$
so
$$
z=L(x-x',0) + (w+w') \in \mathrm{Im}(X\times \{0\}) + W.
$$

\item Claim:  $L(X\times \{0\})+W =Z$ $\Rightarrow$ $\pi$ is surjective.

Given any $y\in Y$, we have $L(0,y)\in Z$, so
$L(0,y)=L(x,0)+w$ for some $x\in X$ and $w\in W$. We have
$L(-x,y)=w$ so $(-x,y)\in L^{-1}(W)$ and $\pi(-x,y)=y$.
\end{enumerate}
\end{proof}

We now assume that $\chi(\Si)<0$. The formula of $d_\mu$ is given
by \cite[(10.7)]{ym}:
$$
d_\mu=\sum_{\alpha(\mu)>0}(\alpha(\mu)+\tilde{g}-1)\geq 0
$$
where $\tilde{g}\geq 2$ is the genus of $\tSi$. 
Note that $d_\mu\geq 2$ if
$d_\mu\neq 0$, so the real codimension of any lower stratum
in $\cA(P)$  is at least two. Since $\cN_0(P)$ is a deformation
retraction of $\cA_{ss}(P)$,  Proposition \ref{thm:path} implies
the following.
\begin{thm}\label{thm:Nconnected}
Let $\Si$ be a closed, compact, connected, nonorientable surface with negative Euler characteristic,
or equivalently, $\Si=\Si^\ell_i$ where $i=1,2$ and $\ell\geq 1$.
Let $G$ be a compact, connected Lie group, and let $P$ be a principal $G$-bundle over $\Si$.
Then the space $\cN_0(P)$ of flat connections on $P$ is nonempty and connected.
\end{thm}

\begin{cor}\label{thm:Mconnected}
Let $P$ be as in Theorem \ref{thm:Nconnected}.
Then the moduli space $\cM(P)$ of gauge equivalence classes of 
flat connections on $P$ is nonempty and connected.
\end{cor}

Note that the connectedness of $\cN_0(P)$ implies the
connectedness of $\cM(P)=\cN_0(P)/\cG(P)$, but not vice versa,
so in general Theorem \ref{thm:Nconnected} is stronger
than Corollary \ref{thm:Mconnected}.

Corollary \ref{thm:Mconnected} extends \cite[Theorem 5.2]{HL3} to the case $\Si=4\RP^2$.
We thank the referee of \cite{HL1} for suggesting this approach to us.
During the revision of this paper, D. Ramras obtained the
following extension of Theorem \ref{thm:Nconnected} in the case $G=U(n)$
(see \cite[Proposition 4.9]{R}):
\begin{thm}\label{thm:connectivity}
Let $\Si$ be a closed, compact, connected, nonorientable surface which is not
$\RP^2$ $($or equivalently, $\chi(\Si)\leq 0)$.
Let $P$ be a principal $U(n)$-bundle $(n\geq 2)$ over $\Si$, and let
$\cN_0(P)$ be the space of flat connections on $P$. Then
$\cN_0(P)$ is $(1-\chi(\Si))(n-1)-1$ connected.
\end{thm}

\section{$U(n)$-Connections on Orientable Surfaces} \label{sec:Uorientable}

\subsection{Connected components of the representation variety and their reductions}
Any point in
$$
\Hom(\Gamma_\bR(\Si^\ell_0), U(n))/U(n)\cong \ymU{\ell}{0}/U(n)
$$
can be represented by $(V,X)\in U(n)^{2\ell}\times \fu(n)$, where
$X$ is a {\em diagonal} matrix. Actually, there is
a unique representative such that
$$
X=-2\pi \sqrt{-1} \left(
\begin{array}{ccc}
\mu_1 & & 0\\
& \ddots & \\
0& &\mu_n
\end{array}
\right)
$$
where $\mu_1\geq \cdots \geq \mu_n$. Suppose that
$$
\mu=(\mu_1,\ldots,\mu_n)= (
\underbrace{\lambda_1,\cdots,\lambda_1}_{n_1},\ldots,
\underbrace{\lambda_m,\ldots,\lambda_m}_{n_m} )
$$
where $\lambda_1 >\cdots >\lambda_m$ and
$n_1+\cdots + n_m= n$. Then
$$
U(n)_X= U(n_1)\times \cdots \times U(n_m),
$$
and
$$
\exp(X)=\fm(V)\in SU(n_1)\times \cdots\times SU(n_m)
$$
where $\fm(\ab)=\pab$ as before.

Given an $n_j\times n_j$ matrix $A_j$, let $\diag(A_1,\ldots,A_m)$ denote
the $n\times n$ matrix
$$
\left(
\begin{array}{ccc}
A_1 & & 0\\
& \ddots & \\
0& &A_m
\end{array}
\right)
$$
where $n=n_1+\cdots+n_m$. With this notation, we have
\begin{eqnarray*}
X&=& -2\pi\sqrt{-1}\diag\left(\lambda_1 I_{n_1},\ldots, \lambda_m I_{n_m}\right),\\
\exp(X)&=&\diag\bigl(e^{-2\pi\sqrt{-1}\lambda_1} I_{n_1},\ldots,e^{-2\pi\sqrt{-1}\lambda_m} I_{n_m}\bigr),
\end{eqnarray*}
where $e^{-2\pi\sqrt{-1}\lambda_j} I_{n_j}\in SU(n_j)$, or equivalently, $e^{-2\pi\sqrt{-1}n_j\lambda_j}=1$.
So
$$
k_j =\lambda_j n_j\in\bZ
$$
and
$$
X=-2\pi\sqrt{-1}\diag\Bigl(\frac{k_1}{n_1} I_{n_1},\ldots,\frac{k_m}{n_m} I_{n_m}\Bigr).
$$

For each pair $(n,k)\in \bZ_{>0}\times \bZ$, define
\begin{eqnarray*}
I_{n,k}&=&\Bigl \{ \mu=
(\mu_1,\ldots,\mu_n)=  \Bigl(
\underbrace{\frac{k_1}{n_1},\ldots,\frac{k_1}{n_1}}_{n_1},\ldots,
\underbrace{\frac{k_m}{n_m},\ldots,\frac{k_m}{n_m}}_{n_m}\Bigr )\Bigr | \\
&& \Bigl. n_j\in \bZ_{>0}, k_j\in \bZ, \sum_{j=1}^m n_j= n, \sum_{j=1}^m k_j=k,
\frac{k_1}{n_1} >\cdots >\frac{k_m}{n_m} \Bigr \}
\end{eqnarray*}
Given
\begin{equation}\label{eqn:mu}
\mu=(\mu_1,\ldots,\mu_n)=  \bigl(
\underbrace{\frac{k_1}{n_1},\ldots,\frac{k_1}{n_1}}_{n_1},\ldots,
\underbrace{\frac{k_m}{n_m},\ldots,\frac{k_m}{n_m}}_{n_m}\bigr )
\in I_{n,k},
\end{equation}
let
\begin{equation}\label{eqn:Xmu}
X_\mu=-2\pi\sqrt{-1}\diag\Bigl(\frac{k_1}{n_1} I_{n_1},\ldots,\frac{k_m}{n_m} I_{n_m}\Bigr),
\end{equation}
and let $C_\mu$ be the conjugacy class of $X_\mu$.

Note that if $(V, X)\in \ymU{\ell}{0}$,
then $X\in C_\mu$ for some $\mu \in \bigcup_{k\in\bZ}I_{n,k}$.
From now on, we identify $\Hom(\Gamma_\bR(\Si^\ell_0),U(n))$ with $\ymU{\ell}{0}$.

Given $\mu\in \bigcup_{k\in\bZ}I_{n,k}$, define
$$
\ymU{\ell}{0}_\mu= \{ (V,X)\in U(n)^{2\ell}\times C_\mu \mid
V\in (U(n)_X)^{2\ell},\ \fm(V)=\exp(X)\bigr\}.
$$
Then
$$
\ymU{\ell}{0}=\bigcup_{k\in\bZ}\bigcup_{\mu\in I_{n,k}}\ymU{\ell}{0}_\mu.
$$
The $G$-action on $\ymU{\ell}{0}$ preserves $\ymU{\ell}{0}_\mu$.
We will show that
\begin{pro}\label{thm:Xzero}
$$
\{ \ymU{\ell}{0}_\mu \mid \mu\in \bigcup_{k\in\bZ} I_{n,k}\}
$$
are the connected components of
$\ymU{\ell}{0}=\Hom(\Gamma_\bR(\Si^{\ell}_0),U(n))$.
Given
\begin{equation}\label{eqn:muII}
\mu=(\mu_1,\ldots,\mu_n)=  \bigl(
\underbrace{\frac{k_1}{n_1},\ldots,\frac{k_1}{n_1}}_{n_1},\ldots,
\underbrace{\frac{k_m}{n_m},\ldots,\frac{k_m}{n_m}}_{n_m}\bigr )
\in I_{n,k}.
\end{equation}
we have a homeomorphism
\begin{equation}\label{eqn:ymUzero}
\ymU{\ell}{0}_\mu/U(n)=\prod_{j=1}^m X_{\mathrm{YM}}^{\ell,0}(U(n_j))_{\frac{k_j}{n_j},\ldots,\frac{k_j}{n_j}}/U(n_j)
\end{equation}
\end{pro}
\begin{proof}
Let $\mu$ be as in \eqref{eqn:muII}.
Let $\pi: \ymU{\ell}{0}\to C_\mu$ be defined
by
$
(V,X)\mapsto X.
$
Then $\pi$ is a fibration, so there is a homeomorphism
\begin{equation}\label{eqn:piX}
\ymU{\ell}{0}_\mu/U(n) = \pi^{-1}(X_\mu)/U(n)_{X_\mu},
\end{equation}
where $U(n)_{X_\mu}=U(n_1)\times\cdots \times U(n_m)$. We have
\begin{eqnarray*}
\pi^{-1}(X_\mu) &\cong&  \{V\in U(n)_{X_\mu}^{2\ell} \mid \fm(V)=\exp(X_\mu) \}\\
&\cong&\prod_{j=1}^m \{ V\in U(n_j)^{2\ell}\mid \fm(V)=e^{-2\pi\sqrt{-1}k_j/n_j} I_{n_j}\}\\
&=& \prod_{j=1}^m X_{\mathrm{YM}}^{\ell,0}(U(n_j))_{\frac{k_j}{n_j},\ldots,\frac{k_j}{n_j}}.
\end{eqnarray*}
This proves \eqref{eqn:ymUzero}. The set
$$
X_{\mathrm{YM}}^{\ell,0}(U(n_j))_{\frac{k_j}{n_j},\ldots,\frac{k_j}{n_j}}
= \{ V\in U(n_j)^{2\ell}\mid \fm(V)=e^{-2\pi\sqrt{-1}k_j/n_j} I_{n_j}\}
$$
is nonempty and connected by \cite[Theorem 3]{HL3}. So $\pi^{-1}(X_\mu)$
is nonempty and connected. Together with \eqref{eqn:piX}, this implies that
$\ymU{\ell}{0}_\mu$ is nonempty and connected.

Define  $T:\ymU{\ell}{0}\to \bR^n$ by
$$
T(V,X)=\left(\frac{\sqrt{-1}}{2\pi}\Tr X,\left(\frac{\sqrt{-1}}{2\pi}\right)^2\Tr(X^2),\ldots,
\left(\frac{\sqrt{-1}}{2\pi}\right)^n\Tr(X^n)\right)
$$
where $V\in U(n)^{2\ell}$ and $X\in \fu(n)$. The characteristic polynomial
$$
P_X(t)=\det(tI-X)=(t+2\pi\sqrt{-1}\mu_1)\cdots (t+2\pi\sqrt{-1}\mu_n)
$$
of a matrix $X$ is determined by $\Tr X,\Tr(X^2),\ldots, \Tr(X^n)$ and the conjugacy
class of $X$ in $\fu(n)$ is determined by $P_X(t)$, so
$T(V,X)=T(V',X')$
if and only if $X$ and $X'$ are in the same conjugacy class.

Given $\mu\in \bigcup_{k\in\bZ}I_{n,k}$,
define
$$
v_\mu=\Bigl(\sum_{i=1}^n \mu_i, \sum_{i=1}^n \mu_i^2,\ldots,\sum_{i=1}^n \mu_i^n\Bigr)\in \bQ^n.
$$
Note that $v_{\mu}=v_{\mu'}$ if and only if $\mu=\mu'$.

The function $T$ is a continuous function, and its image
$\{ v_\mu\mid \mu\in \bigcup_{k\in \bZ} I_{n,k} \}$ is a discrete set, so
$$
\{ \ymU{\ell}{0}_\mu = T^{-1}(v_\mu) \mid \mu\in \bigcup_{k\in\bZ}I_{n,k} \}
$$
are connected components of $\ymU{\ell}{0}$.
\end{proof}

\subsection{Equivariant Morse stratification and equivariant Poincar\'{e} series}
Let $P^{n,k}$ be the topological principal $U(n)$-bundle on $\Si^\ell_0$ with
$c_1(P^{n,k})=k[\omega]\in H^2(\Si^\ell_0)$. Let $\cA^{n,k}$ be the space of $U(n)$-connections
on $P^{n,k}$, and let $\cN^{n,k}\subset \cA^{n,k}$ be the space of Yang-Mills
$U(n)$-connections on $P^{n,k}$. Let $\cG^{n,k}$ be the group of gauge transformations
on $P^{n,k}$, and let $\cG_0^{n,k}$ be the subgroup of gauge transformations
which take value the identity $e$ at a fixed point $x_0\in\Si^\ell_0$.
We have
$$
\Hom(\Gamma_\bR(\Si^\ell_0),U(n))=\bigcup_{k\in\bZ}\Hom(\Gamma_\bR(\Si^\ell_0),U(n))_k,
$$
where $\Hom(\Gamma_\bR(\Si^\ell_0),U(n))_k \cong  \cN^{n,k}/ \cG^{n,k}_0$.
The connected components of $\Hom(\Gamma_\bR(\Si^\ell_0), U(n))_k$ are
$$
\{ \ymU{\ell}{0}_\mu\mid \mu\in I_{n,k} \}.
$$
Let $\cN_\mu$ be the preimage of $\ymU{\ell}{0}_\mu$ under the projection
$$
\cN^{n,k}\to \cN^{n,k}/\cG^{n,k}_0\cong \Hom(\Gamma_\bR(\Si^\ell_0), U(n))_k.
$$
so that $\ymU{\ell}{0}_\mu =\cN_\mu/\cG_0(P)$.

We fix $(n,k)\in \bZ_{>0}\times \bZ$, and write $\cG=\cG^{n,k}$.
Let $\cA_\mu$ be the stable manifold of $\cN_\mu$ with respect to
Yang-Mills functional. Then
$$
\cA^{n,k}=\bigcup_{\mu\in I_{n,k}} \cA_\mu.
$$
is the $\cG$-equivariant Morse stratification
of $\cA^{n,k}$ given by the Yang-Mills functional \cite{ym,da,Ra,CW}.
Let $\lambda_\mu$ be the real codimension
of $\cA_\mu$ in $\cA^{n,k}$. It was computed in \cite[Section 7]{ym} that
$$
\lambda_\mu=2 d_\mu,\ \
d_\mu=\sum_{i<j}(\mu_i-\mu_j+(g-1)).
$$
The gradient flow of the Yang-Mills functional gives a $\cG$-equivariant
deformation retraction $\cA_\mu \to \cN_\mu$. For the purpose of equivariant cohomology,
the following equivariant pairs are equivalent:
$$
(\cA_\mu,\cG)\sim (\cN_\mu,\cG)\sim (\ymU{\ell}{0}_\mu, U(n)).
$$
In other words, we have the following homotopy equivalences of homotopic orbit spaces:
$$
{\cA_\mu}^{h \cG}\sim {\cN_\mu}^{h \cG}\sim {\ymU{\ell}{0}_\mu}^{h U(n)}.
$$
Together with the reduction Proposition \ref{thm:Xzero},
we conclude that
\begin{thm}
Let $K$ be a field. For any
$$
\mu=(\mu_1,\ldots,\mu_n)=  \bigl(
\underbrace{\frac{k_1}{n_1},\ldots,\frac{k_1}{n_1}}_{n_1},\ldots,
\underbrace{\frac{k_m}{n_m},\ldots,\frac{k_m}{n_m}}_{n_m}\bigr )
\in I_{n,k},
$$
we have
\begin{equation}\label{eqn:Hzero}
H^*_{\cG}(\cA_\mu;K)=H^*_{\cG}(\cN_\mu;K)=
H^*_{U(n)}(\ymU{\ell}{0}_\mu;K) \cong \bigotimes_{j=1}^m H^*_{U(n_j)}(V_{ss}(P^{n_j,k_j});K),
\end{equation}
\begin{equation}\label{eqn:Pzero}
 P^\cG_t(\cA_\mu;K)=P^\cG_t(\cN_\mu;K)=
P^{U(n)}_t(\ymU{\ell}{0}_\mu;K)=\prod_{i=1}^m P^{U(n_j)}_t(V_{ss}(P^{n_j,k_j});K).
\end{equation}
\end{thm}

\subsection{Involution}\label{sec:ZU}

Given $\mu\in I_{n,k}$, define
\begin{eqnarray*}
\ZymU{1}_\mu&=&\{(V,c,\bV,\bc,X)\in U(n)^{2(2\ell+1)}\times C_\mu \mid
V , c\bV c^{-1}\in (U(n)_X)^{2\ell},\\
&& \quad \fm(V)=\exp(X/2)c\bc,\  \fm(\bV)=\bc \exp(-X/2)c\}\\
\ZymU{2}_\mu&=&\{(V,d,c,\bV,\bd,\bc,X)\in U(n)^{2(2\ell+2)}\times C_\mu\mid
V, d^{-1}c\bV c^{-1}d \in (U(n)_X)^{2\ell},\\
&& \quad d^{-1},c\bc \in U(n)_X,\
\fm(V)=\exp(X/2)c\bd c^{-1}d,\  \fm(\bV)=\bc d\exp(-X/2) \bc^{-1}\bd \}
\end{eqnarray*}

Then for $i=1,2$,
$$
\ZymU{i}=\bigcup_{k\in\bZ}\bigcup_{\mu\in I_{n,k}} \ZymU{i}_\mu
$$
$$
\Phi^{\ell,i}_{U(n)}(\ZymU{i}_\mu)=\ymU{2\ell+i-1}{0}_\mu.
$$

Define $\tau_0:I_{n,k}\to I_{n,-k}$ by
$$
(\mu_1,\mu_2,\ldots,\mu_n)\mapsto (-\mu_n,\cdots,-\mu_2,-\mu_1).
$$
It is easy to check that if $X\in C_\mu$ then $-\Ad(\bc)(X)\in C_{\tau_0(\mu)}$. So
$$
\tau(\ZymU{i}_\mu)= \ZymU{i}_{\tau_0(\mu)}.
$$
Thus we conclude that:
\begin{thm}
The set
$$
\ZymU{i}^\tau_\mu=\ZymU{i}^{\tau} \cap \ZymU{i}_\mu
$$
is nonempty if and only if $X_\mu$ is conjugate to $-X_\mu$, i.e. $\tau_0(\mu)=\mu$.
In other words, if we define
$$
I_n=I_{n,0}^{\tau_0}=\{\mu\in I_{n,0}\mid \tau_0(\mu)=\mu\}.
$$
Then
$$
\ZymU{i}^\tau=\bigcup_{\mu\in I_n}
\ZymU{i}_\mu^{\tau}.
$$
\end{thm}

\section{$U(n)$-Connections on Nonorientable Surfaces}\label{sec:Unonorientable}

\subsection{Connected components of the representation variety and their reductions}

Given $\mu\in I_{n,k}$, let $C_{\mu/2}$ denote the conjugacy class of $X_\mu/2$
in $\fu(n)$.
$\ZymU{1}_\mu^\tau$ can be identified with
\begin{eqnarray*}
\ymU{\ell}{1}_\mu&=& \{ (V,c,X)\in U(n)^{2\ell+1}\times C_{\mu/2}\mid
 V\in (U(n)_X)^{2\ell},
 \\&&\quad
 \Ad(c)(X)=-X,\ \fm(V) =\exp(X)c^2\}.
\end{eqnarray*}
while $\ZymU{2}_\mu^\tau$ can be identified with
\begin{eqnarray*}
\ymU{\ell}{2}_\mu&=& \{ (V,d,c,X)\in U(n)^{2\ell+2}\times C_{\mu/2}\mid
V\in (U(n)_X)^{2\ell},\ d\in U(n)_X,\\
&&\quad  \Ad(c)(X)=-X,\ \fm(V)=\exp(X) cd c^{-1} d\}.
\end{eqnarray*}

Note that a Yang-Mills connection on a principal
$U(n)$-bundle $P$ over a nonorientable surface $\Si$
induces a flat connection on the $U(1)$-bundle $\det(P)$.
More explicitly, define
$$
\det: U(n)^{2\ell}\to U(1)^{2\ell},\quad
(\ab)\mapsto (\det(a_1),\det(b_1),\ldots,\det(a_\ell),\det(b_\ell))\in U(1)^{2\ell}.
$$
We have
\begin{eqnarray*}
\det: \ymU{\ell}{1}\to \flS{\ell}{1} &&
(V,c,X)\mapsto (\det(V),\det(c))\\
\det: \ymU{\ell}{2}\to \flS{\ell}{2} &&
(V,d,c,X)\mapsto (\det(V),\det(d),\det(c))
\end{eqnarray*}
where $\flS{\ell}{1}$ and $\flS{\ell}{2}$ are as in Example \ref{flS}.
Let
\begin{eqnarray*}
\flS{\ell}{1}^{\pm 1}&=&\{(\ab,c)\in U(1)^{2\ell+1} \mid c=\pm1\}\\
\flS{\ell}{2}^{\pm 1}&=&\{(\ab,d,c)\in U(1)^{2\ell+2} \mid d=\pm 1 \}
\end{eqnarray*}
Then $\flS{\ell}{i}^{+1}$ and $\flS{\ell}{i}^{-1}$ are the two
connected components of $\flS{\ell}{i}$. Let
\begin{eqnarray*}
\ymU{\ell}{i}^{\pm 1} &=& (\det)^{-1}\left(\flS{\ell}{i}^\pm\right),\\
\ymU{\ell}{i}_\mu^{\pm 1} &=&\ymU{\ell}{i}_\mu \cap \ymU{\ell}{i}^\pm.
\end{eqnarray*}
Then
$$
\ymU{\ell}{i}_\mu =\ymU{\ell}{i}_\mu^{+1} \cup \ymU{\ell}{i}_\mu^{-1}.
$$

Any $\mu\in I_n$ is of the following form:
$$
\mu = (\nu,\underbrace{0,\ldots,0}_{n_0},\tau_0(\nu))
$$
where $\nu \in I_{n',k}$, $\tau_0(\nu)\in I_{n',-k}$,
$n'\geq 0$, $n_0\geq 0$, $k>0$.

We will show that
\begin{pro}\label{thm:Ncomponents}
Let $\mu=(\nu,\underbrace{0,\ldots,0}_{n_0},\tau_0(\nu))\in I_n=I^{\tau_0}_{n,0}$, where
$$\nu\in I_{n',k},\quad  n', n_0\geq 0,\quad  2n'+n_0=n,\quad  k>0.$$
\begin{enumerate}
\item[(i)] Suppose that $n_0>0$. For $i=1,2$,
$\ymU{\ell}{i}_\mu^{+1}$, $\ymU{\ell}{i}_\mu^{-1}$ are nonempty and connected
for $\ell\geq 1$.
\item[(ii)]Suppose that $n_0=0$ so that $2n'=n>0$. For $i=1,2$,
\begin{equation}\label{eqn:pm}
\ymU{\ell}{i}_\mu =
\ymU{\ell}{i}_\mu^{(-1)^{n'i+k}}
\end{equation}
$\ymU{\ell}{i}_\mu$ is nonempty and connected unless $i=1$ and $\ell=0$.
\end{enumerate}
\end{pro}

\begin{proof}
(i) $n_0>0$.
$\nu$ is of the form
$$
\nu=\Bigl(\underbrace{\frac{k_1}{n_1},\ldots \frac{k_1}{n_1}}_{n_1},
\ldots,
\underbrace{\frac{k_m}{n_m},\ldots,\frac{k_m}{n_m}} _{n_m}
\Bigr)
$$
where
$$
\frac{k_1}{n_1}>\ldots>\frac{k_m}{n_m}>0,
$$
and
$$
X_\mu=-2\pi\sqrt{-1}\diag\Bigl(\frac{k_1}{n_1}I_{n_1},\ldots,\frac{k_m}{n_m}I_{n_m}, 0 I_{n_0},
-\frac{k_m}{n_m}I_{n_m},\ldots, -\frac{k_1}{n_1}I_{n_1}\Bigr).
$$

Let $\pi_1:\ymU{\ell}{1}_\mu \to C_{\mu/2}$ be defined by
$(V,c,X)\mapsto X$,
and let  $\pi_2:\ymU{\ell}{2}_\mu \to C_{\mu/2}$ be defined by
$(V,d,c,X)\mapsto X$.
Then $\pi_1$ and $\pi_2$ are fibrations. So
$$
\ymU{\ell}{i}_\mu/U(n)\cong \pi_i^{-1}(X_\mu/2)/U(n)_\mu
$$
where
$$
U(n)_\mu = U(n)_{X_\mu/2}=U(n)_{X_\mu}.
$$

Let
$$
e_\mu=\left(\begin{array}{ccccccc}
& & & & & & I_{n_1} \\
&0 & & & & \cdot & \\
& & & & I_{n_m} & & \\
& & & I_{n_0} & & &\\
& & I_{n_m} & & & &\\
& \cdot & & & &0 & \\
I_{n_1} & & & & & &
 \end{array}\right)\in U(n)
$$
Then
$$
\Ad(e_\mu)X_\mu =-X_\mu,\quad \det(e_\mu)=(-1)^{n'},\quad e_\mu^2=I_n.
$$

We have $\pi_i^{-1}(X_\mu/2)\cong V^i_\mu$, where
\begin{eqnarray*}
V_\mu^1&=& \{ (V,c')\in U(n)_\mu^{2\ell+1}\mid
\fm(V) =\exp(X_\mu/2) e_\mu c' e_\mu c'\},\\
V_\mu^2&=& \{ (V,d,c')\in U(n)_\mu^{2\ell+2}\mid
\fm(V) =\exp(X_\mu/2) e_\mu c'd (e_\mu c')^{-1}d\}.
\end{eqnarray*}
Under the identification $\pi_i^{-1}(X_\mu/2)\cong V^i_\mu$, the actions of $U(n)_{\mu}$
on $V^1_\mu$ and $V^2_\mu$ are given by
$$
g\cdot(V,c')= (gVg^{-1}, (e_\mu g e_\mu) c' g^{-1}) \quad \textup{and}\quad
g\cdot(V,d,c')=(gVg^{-1}, gdg^{-1},  (e_\mu g e_\mu) c' g^{-1})
$$
respectively, where $V\in U(n)_\mu^{2\ell}$, and $g,c',d\in U(n)_\mu$.

Any $a_i,b_i,c',d\in U(n)_\mu$ are of the form
\begin{eqnarray*}
a_i&=&\diag(A_1^i,\ldots,A_m^i,A^i,\bar{A}_m^i,\ldots,\bar{A}_1^i)\\
b_i&=&\diag(B_1^i,\ldots,B_m^i,B^i,\bar{B}_m^i,\ldots,\bar{B}_1^i)\\
c'&=&\diag(\bar{C}_1,\ldots\bar{C}_m,C,C_m,\ldots,C_1)\\
d &=& \diag(D_1,\ldots D_m,D,\bar{D}_m,\ldots,\bar{D}_1)
\end{eqnarray*}
where
$$
A^i_j,\ \bar{A}^i_j,\  B^i_j,\ \bar{B}^i_j,\ C_j,\ \bar{C}_j,\  D_j,\ \bar{D}_j\in U(n_j),
\quad \quad
A^i,\ B^i,\ C,\ D\in U(n_0).
$$

If $(V,c')\in V^1_\mu$, we have
\begin{equation}\label{eqn:C}
\prod_{i=1}^\ell[A^i_j,B^i_j]=e^{\frac{-\pi\sqrt{-1}k_j}{n_j}} C_j \bar{C}_j,\quad
\prod_{i=1}^\ell[\bar{A}^i_j,\bar{B}^i_j]= e^{\frac{\pi\sqrt{-1}k_j}{n_j}} \bar{C}_j C_j
\end{equation}
for $j=1,\ldots,m$, and
\begin{equation}\label{eqn:Cflat}
\prod_{i=1}^\ell[A^i,B^i]=C^2.
\end{equation}
By \eqref{eqn:C}, we have
$$
1=\det\Bigl(e^{\frac{-\pi\sqrt{-1}k_j}{n_j} } C_j\bar{C}_j\Bigr)=(-1)^{k_j}\det(C_j)\det(\bar{C}_j).
$$
By \eqref{eqn:Cflat}, we have $1=\det(C^2)=\det(C)^2$.
 Recall that $\det(c)=\det(e_\mu)\det(c)=(-1)^{n'}\det(c')$, so
\begin{equation}\label{eqn:detC}
\det(c)=(-1)^{n'+k}\det(C), \quad \det(C)=\pm 1.
\end{equation}
Note that the equations in \eqref{eqn:C} are exactly the defining equations of
$Z_{\mathrm{YM}}^{\ell,1}(U(n_j))_{\frac{k_j}{n_j},\ldots,\frac{k_j}{n_j}}$, where
$\ZymU{1}_\mu$  is defined as in Section \ref{sec:ZU}, and   \eqref{eqn:Cflat}
is exactly the defining equation for $X_{\mathrm{flat}}^{\ell,1}(U(n_0))$.
So we have the following homeomorphism:
\begin{equation}\label{eqn:Vone}
V^1_\mu\cong  X^{\ell,1}_{\mathrm{flat}}(U(n_0)) \times \prod_{j=1}^m
Z^{\ell,1}_{\mathrm{YM}}(U(n_j))_{\frac{k_j}{n_j},\ldots,\frac{k_j}{n_j}}.
\end{equation}

If $(V,d,c')\in V^2_\mu$, we have
\begin{equation}\label{eqn:CD}
\prod_{i=1}^\ell[A^i_j,B^i_j]=e^{\frac{-\pi\sqrt{-1}k_j}{n_j}} C_j \bar{D}_j C_j^{-1} D_j,\quad
\prod_{i=1}^\ell[\bar{A}^i_j,\bar{B}^i_j]= e^{\frac{\pi\sqrt{-1}k_j}{n_j}}
\bar{C}_j D_j\bar{C}_j^{-1} \bar{D}_j
\end{equation}
for $j=1,\ldots,m$, and
\begin{equation}\label{eqn:CDflat}
\prod_{i=1}^\ell[A^i,B^i]=C D C^{-1}D.
\end{equation}
By \eqref{eqn:CD}, we have
$$
1=\det\Bigl(e^{\frac{-\pi\sqrt{-1}k_j}{n_j} } D_j\bar{D}_j\Bigr)=(-1)^{k_j}\det(D_j)\det(\bar{D}_j).
$$
By \eqref{eqn:CDflat}, we have $1=\det(D^2)=\det(D)^2$. We conclude that
\begin{equation}\label{eqn:detD}
\det(d)=(-1)^k \det(D),\quad \det(D)=\pm 1.
\end{equation}
Note that the equations in \eqref{eqn:CD} are exactly the defining equations of
$Z_{\mathrm{YM}}^{\ell,2}(U(n_j))_{\frac{k_j}{n_j},\ldots,\frac{k_j}{n_j}}$,
where $\ZymU{2}_\mu$  is defined as in Section \ref{sec:ZU}, and
\eqref{eqn:CDflat} is exactly the defining equation of
$X_{\mathrm{flat}}^{\ell,2}(U(n_0))$. So we have the following homeomorphism:
\begin{equation}\label{eqn:Vtwo}
V^2_\mu\cong  X^{\ell,2}_{\mathrm{flat}}(U(n_0)) \times \prod_{j=1}^m
Z^{\ell,2}_{\mathrm{YM}}(U(n_j))_{\frac{k_j}{n_j},\ldots,\frac{k_j}{n_j}}.
\end{equation}

The homeomorphisms \eqref{eqn:Vone}, \eqref{eqn:Vtwo} are $U(n)_\mu$-equivariant: for $i=1,2$,
it is straightforward to check that the action of
$$
U(n)_\mu= U(n_1)\times \ldots\times U(n_m)\times U(n_0)\times U(n_m)\times \ldots\times U(n_1)
$$
on $V^i_\mu$ is compatible with the
action of $U(n_0)$ on   $X^{\ell,i}_{\mathrm{flat}}(U(n_0))$ and the actions of $U(n_j)^2$ on
$Z^{\ell,i}_{\mathrm{YM}}(U(n_j))_{\frac{k_j}{n_j},\ldots,\frac{k_j}{n_j}}$. So
\begin{eqnarray*}
\ymU{\ell}{i}_\mu /U(n)&\cong& V_\mu^i/U(n)_\mu\\
&\cong & X^{\ell,i}_{\mathrm{flat}}(U(n_0))/U(n_0)\times
\prod_{j=1}^m \Bigl( Z^{\ell,i}_{\mathrm{YM}}(U(n_j))_{\frac{k_j}{n_j},\ldots,\frac{k_j}{n_j}}/U(n_j)^2 \Bigr) \\
&\cong& X^{\ell,i}_{\mathrm{flat}}(U(n_0))/U(n_0)\times
\prod_{j=1}^m \Bigl( X^{2\ell+i-1,0}_{\mathrm{YM}}(U(n_j))_{\frac{k_j}{n_j},\ldots,\frac{k_j}{n_j}}/U(n_j)\Bigr)\\
&\cong& X^{\ell,i}_{\mathrm{flat}}(U(n_0))/U(n_0)\times
 X^{2\ell+i-1,0}_{\mathrm{YM}}(U(n'))_\nu/U(n').
\end{eqnarray*}
From \eqref{eqn:detC} and \eqref{eqn:detD}, we see that
$$
\ymU{\ell}{i}_\mu^\pm 1 /U(n)\cong
X_{\mathrm{flat}}^{\ell,i}(U(n_0))^{\pm(-1)^{n'i+k} }\times  X^{2\ell+i-1,0}_{\mathrm{YM}}(U(n'))_\nu/U(n')
$$
where
$$
X_{\mathrm{flat}}^{\ell,i}(U(n_0))^{\pm 1} = X_{\mathrm{YM}}^{\ell,i}(U(n_0))_{0,\ldots,0}^{\pm 1}.
$$
Recall from Proposition \ref{thm:Xzero} that if $2\ell+i-1\geq 1$ then
$$
X^{2\ell+i-1,0}_{\mathrm{YM}}(U(n'))_\nu
$$
is nonempty and connected. By Theorem \ref{thm:components},
$$
X_{\mathrm{flat}}^{\ell,i}(U(n_0))^{+1}, \quad X_{\mathrm{flat}}^{\ell,i}(U(n_0))^{-1}
$$
are nonempty and connected for $\ell\geq 1$. We conclude that
$\ymU{\ell}{i}_\mu^{+1}$ and $\ymU{\ell}{i}_\mu^{-1}$ are nonempty and connected for
$\ell\geq 1$.

\noindent
(ii) $n_0=0$. The calculations in this case are the same as those in (i), except that we do not
have the factor $U(n_0)$ so we do not have the matrices $A^i,B^i,$$C,D$. We conclude that
$\det(c)=(-1)^{n'+k}$ when $i=1$ and  $\det(d)=(-1)^k$ when $i=2$. So
$$
\ymU{\ell}{i}_\mu =
\ymU{\ell}{i}_\mu^{(-1)^{n'i+k}}
$$
and
$$
\ymU{\ell}{i}_\mu/U(n)\cong X^{2\ell+i-1,0}_{\mathrm{YM}}(U(n'))_\nu/U(n').
$$
Recall that $X^{2\ell+i-1,0}_{\mathrm{YM}}(U(n'))_\nu$ is nonempty and
connected if $2\ell+i-1\geq 1$. So for $i=1,2$,
$\ymU{\ell}{i}_\mu$ is nonempty and connected unless $i=1$ and $\ell=0$.
\end{proof}

Let $I_n=I_{n,0}^{\tau_0}$. Then $I_n=I_n^0\cup I_n^{i,+}\cup I_n^{i,-}$, where
\begin{eqnarray*}
I_n^0&=& \{\mu\in I_n, \mu_i=0 \textup{ for some }i \},\\
I_n^{i,+} &=& \{ \mu=(\nu,\tau_0(\nu))\mid \nu\in I_{n',k}, n'i+k \textup{ is even} \},\\
I_n^{i,-} &=& \{ \mu=(\nu,\tau_0(\nu))\mid \nu\in I_{n',k}, n'i+k \textup{ is even} \},
\end{eqnarray*}
where $i=1,2$.
When $n$ is odd, we have $I_n=I_n^0$.

The proof of Proposition \ref{thm:Ncomponents} gives the following.
\begin{pro}\label{thm:Xcomponents}
Suppose that  $i=1,2$ and $\ell\geq 1$.
\begin{itemize}
\item[(i)]If $n$ is odd, the connected components of
$\ymU{\ell}{i}=\Hom(\Gamma_\bR(\Si^\ell_i),U(n))$ are
$$
\{ \ymU{\ell}{i}_\mu^{+1}  \mid \mu\in I_n^0 \}\cup
\{ \ymU{\ell}{i}_\mu^{-1} \mid \mu \in I_n^0 \}.
$$
\item[(ii)]If $n$ is even, the connected components of
$\ymU{\ell}{i}=\Hom(\Gamma_\bR(\Si^\ell_i),U(n))$ are
\begin{eqnarray*}
 && \quad \quad \quad \{ \ymU{\ell}{i}_\mu^{+1} \mid \mu\in I_n^0\}\cup
\{\ymU{\ell}{i}_\mu=\ymU{\ell}{i}_\mu^{+1}\mid \mu\in I_n^{i,+}\} \\
&& \cup
\{\ymU{\ell}{i}_\mu^{-1} \mid  \mu\in I_n^0\}\cup
\{\ymU{\ell}{i}_\mu=\ymU{\ell}{i}_\mu^{-1}\mid \mu\in I_n^{i,-}\}.
\end{eqnarray*}
\end{itemize}
\end{pro}

\begin{pro}\label{thm:Xreduce_non}
Given
$\mu=(\nu,\underbrace{0,\ldots,0}_{n_0},\ \tau_0(\nu))\in I_n$, where
$\nu\in I_{n',k},\  n', n_0\geq 0,\  2n'+n_0=n$,
$$
\nu=\Bigl(\underbrace{\frac{k_1}{n_1},\dots,\frac{k_1}{n_1}}_{n_1},\ldots,
\underbrace{\frac{k_m}{n_m},\ldots,\frac{k_m}{n_m}}_{n_m}\Bigr),\quad
\frac{k_1}{n_1}>\cdots > \frac{k_m}{n_m}>0,
$$
we have the following homeomorphisms.
\begin{enumerate}
\item[(i)]If $n_0>0$, then for $i=1,2$
$$
\ymU{\ell}{i}_\mu^\pm /U(n)
\cong
X_{\mathrm{flat}}^{\ell,i}(U(n_0))^{\pm(-1)^{ni'+k}}/U(n_0)\times \prod_{j=1}^m
\Bigl(X^{2\ell+i-1,0}_{\mathrm{YM}}(U(n_j))_{\frac{k_j}{n_j},\ldots,\frac{k_j}{n_j}}/U(n_j)\Bigr)
$$
\item[(ii)]If $n_0=0$, then for $i=1,2$
$$
\ymU{\ell}{i}_\mu/U(n)=\prod_{j=1}^m
\Bigl(X^{2\ell+i-1,0}_{\mathrm{YM}}(U(n_j))_{\frac{k_j}{n_j},\ldots,\frac{k_j}{n_j}}/U(n_j)\Bigr).
$$
\end{enumerate}
\end{pro}

\subsection{Equivariant Morse stratification and equivariant Poincar\'{e} series}

For $i=1,2$, let $P^{n,+}$ and $P^{n,-}$ be the principal $U(n)$-bundles on
$\Si^\ell_i$ with $c_1(P^{n,+})=0$ and $c_1(P^{n,-})=1$ in
$H^2(\Si^\ell_i;\bZ/2\bZ)\cong \bZ/2\bZ$, respectively.
Let $\cA(\Si^\ell_i)^{n,\pm}$ be the space of connections on $P^{n,\pm}\to \Si^\ell_i$, and let
$\cN(\Si^\ell_i)^{n,\pm}$ denote the space of Yang-Mills $U(n)$-connections
on $P^{n,\pm}\to \Si^\ell_i$. Let $\cG^{n,\pm}$ and $\cG_0^{n,\pm}$ denote the gauge group
and based gauge group, respectively.
In particular,  $\cN(\Si^\ell_i)^{1,\pm}=\cN(\Sigma^\ell_i)_0^{1,\pm}$
is the space of flat $U(1)$-connections
on $P^{1,\pm}\to \Si^\ell_i$, and we have (see \cite{HL2,HL3})
$$
\flS{\ell}{i}^{\pm 1} = \cN(\Si^\ell_i)_0^{1,\pm}/\cG_0^{1,\pm},\quad
\flS{\ell}{i}^{\pm 1}/U(1) =\cN(\Si^\ell_i)_0^{1,\pm} /\cG^{1,\pm}.
$$
Since $c_1(P)= c_1(\det(P))\in H^2(\Si;\bZ/2\bZ)$, we have
$$
\Hom(\Gamma_\bR(\Sigma^\ell_i), U(n))=\ymU{\ell}{i}^{+1}\cup \ymU{\ell}{i}^{-1}
$$
$$
\ymU{\ell}{i}^{\pm 1} = \cN(\Si^\ell_i)^{n,\pm}/\cG_0^{n,\pm},\quad
\ymU{\ell}{i}^{\pm 1}/U(n) = \cN(\Si^\ell_i)^{n,\pm}/\cG^{n,\pm}.
$$

The connected components of $\ymU{\ell}{i}^\pm$ are
$$
\{ \ymU{\ell}{i}_\mu^{\pm 1} \mid \mu\in I_n^0 \}\cup \{\ymU{\ell}{i}_\mu\mid \mu\in I_n^{i,\pm} \}.
$$
When $n$ is odd, both $I_n^{i,+}$ and $I_n^{i,-}$ are empty.
Given $\mu\in I_n^0$, let $\cN(\Si^\ell_i)_\mu^\pm$ be the preimage of
$\ymU{\ell}{i}^{\pm 1}_\mu$ under the projection
$$
\cN(\Si^\ell_i)^{n,\pm}\to \cN(\Si^\ell_i)^{n,\pm}/\cG_0^{n,\pm} = \ymU{\ell}{i}^{\pm 1}.
$$
Given $\mu\in I_n^{i,\pm}$, let $\cN(\Si^\ell_i)_\mu$ be the preimage of
$\ymU{\ell}{i}_\mu$ under the projection
$$
\cN(\Si^\ell_i)^{n,\pm}\to \cN(\Si^\ell_i)^{n,\pm}/\cG_0^{n,\pm} = \ymU{\ell}{i}^{\pm 1}.
$$
In other words,
$$
\cN(\Si^\ell_i)_\mu^\pm/\cG_0^{n,\pm}\cong \ymU{\ell}{i}_\mu^{\pm 1},\ \mu\in I_n^0;\quad
\cN(\Si^\ell_i)_\mu/\cG_0^{n,\pm}\cong \ymU{\ell}{i}_\mu,\ \mu\in I_n^{i,\pm}.
$$
In particular,
$$
\cN(\Si^\ell_i)_{ \tiny{\underbrace{0,\ldots,0}_n} }^{\pm} =\cN(\Si^\ell_i)^{n,\pm}_0
$$
is the space of flat connections on $P^{n,\pm}\to\Si^\ell_i$. The Yang-Mills functional
achieves its absolute minimum $0$ on $\cN(\Si^\ell_i)_0^{n,\pm}$.
The moduli space of gauge equivalence classes of flat connections on
$P^{n,\pm}\to \Si^\ell_i$ is given by
$$
\cM(\Si^\ell_i,P^{n,\pm})=\cN(\Si^\ell_i)_0^{n,\pm}/\cG^{n,\pm}\cong V_{ss}(\Si^\ell_i,P^{n,\pm})/U(n),
$$
where
$$
V_{ss}(\Si^\ell_i,P^{n,\pm})=\flU{\ell}{i}^\pm=\cN(\Si^\ell_i)^{n,\pm}_0/\cG_0^{n,\pm}.
$$

Now assume $\chi(\Si)=2-2\ell-i<0$.
Let $\pi:\tSi\to \Si$ be the orientable double cover.
Then $\pi^*P^{n,\pm}=P^{n,0}\cong \tSi\times U(n)$.
There are involutions $\tau^\pm:P^{n,0}\to P^{n,0}$ which cover the deck transformation
$\tau:\tSi\to \tSi$ such that $P^{n,0}/\tau^\pm = P^{n,\pm}$. Let
$\cA(\Si)^{n,\pm}$ denote the space of connections on $P^{n,\pm}\to \Si$, and let
$\cA(\tSi)^{n,0}$ denote the space of connections on $P^{n,0}\to \tSi$.
Then
$$
\cA(\Si)^{n,\pm}\cong (\cA(\tSi)^{n,0})^{\tau^\pm}.
$$
Given $\mu\in I_n^0$ (resp. $\mu\in I_n^{i,+}\cup I_n^{i,-}$),
let $\cA(\Si)^\pm_\mu$ (resp. $\cA(\Si)_\mu$)be the stable manifold of the critical set
$\cN(\Si)^\pm_\mu$ (resp. $\cN(\Si)_\mu$) of the Yang-Mills functional on $\cA(\Si)^\pm$.
Then for $i=1,2$,
$$
\cA(\Si^\ell_i)^{n,\pm}=\bigcup_{\mu\in I_n^0}\cA(\Si^\ell_i)^\pm_\mu \cup
\bigcup_{\mu\in I_n^{i,\pm}}\cA(\Si^\ell_i)_\mu
$$
is the $\cG$-equivariant Morse stratification of $\cA(\Si^\ell_i)^{n,\pm}$
given by the Yang-Mills functional.

We have
$$
\cN(\Si^\ell_i)^\pm_\mu \subset (\cA(\tSi)^{n,0}_\mu)^{\tau^\pm} \textup{ for } \mu\in I_n^0,\quad
\cN(\Si^\ell_i)_\mu \subset (\cA(\tSi)^{n,0}_\mu)^{\tau^\pm} \textup{ for }\mu\in I_n^{i,\pm}.
$$
By results in \cite{ym} and \cite{da}, $\cA_\mu$ is a complex submanifold of $\cA$. We also know that
$\tau^\pm$ induces anti-holomorphic involution on $\cA$. By Proposition \ref{thm:J}, for
$\mu\in I^0_n$ we have
$$
\mathrm{codim}_\bR\left(\cA(\Si)^\pm_\mu, \cA(\Si)^{n,\pm}\right)=
\mathrm{codim}_\bC\left(\cA(\tSi)_\mu, \cA(\tSi)^{n,0}\right)=\sum_{\alpha<\beta}(\mu_\alpha-\mu_\beta+\chi(\Si^\ell_i)).
$$
The above formula also holds for $\mu\in I_n^{i,\pm}$.

We have the following equivalent equivariant pairs for the purpose
of equivariant (singular) cohomology:
\begin{eqnarray*}
&& (\cA(\Si^\ell_i)^\pm_\mu,\cG^{n,\pm})\sim (\cN(\Si^\ell_i)^\pm_\mu,\cG^{n,\pm})\sim
(\ymU{\ell}{i}^{\pm 1} _\mu,U(n)),\quad \mu \in I_n^0;\\
&& (\cA(\Si^\ell_i)_\mu,\cG^{n,\pm})\sim (\cN(\Si^\ell_i)_\mu,\cG^{n,\pm})\sim
(\ymU{\ell}{i}_\mu,U(n)),\quad \mu \in I_n^{i,\pm}.
\end{eqnarray*}
In other words, we have the following homotopy equivalences
of homotopic orbit spaces:
\begin{eqnarray*}
&& \Bigl(\cA(\Si^\ell_i)^\pm_\mu\Bigr)^{h \cG^{n,\pm}}\sim
\Bigl(\cN(\Si^\ell_i)^\pm_\mu\Bigr)^{h \cG^{n,\pm}}\sim
\Bigl(\ymU{\ell}{i}^{\pm 1}_\mu\Bigr)^{h U(n)},\quad \mu \in I_n^0;\\
&& {\cA(\Si^\ell_i)_\mu}^{h \cG^{n,\pm}}\sim {\cN(\Si^\ell_i)_\mu}^{h \cG^{n,\pm}}\sim
{\ymU{\ell}{i}_\mu}^{h U(n)},\quad \mu \in I_n^{i,\pm}.
\end{eqnarray*}
Together with the reduction Proposition \ref{thm:Xreduce_non}, we have
\begin{thm}
Let $i=1,2$, and let $K$ be a field. Given
$\mu=(\nu,\underbrace{0,\ldots,0}_{n_0},\ \tau_0(\nu))\in I_n$, where
$$
\nu\in I_{n',k},\quad  n', n_0\geq 0,\quad  2n'+n_0=n,
$$
$$
\nu=\Bigl(\underbrace{\frac{k_1}{n_1},\dots,\frac{k_1}{n_1}}_{n_1},\ldots,
\underbrace{\frac{k_m}{n_m},\ldots,\frac{k_m}{n_m}}_{n_m}\Bigr),\quad
\frac{k_1}{n_1}>\cdots > \frac{k_m}{n_m}>0,
$$
we have the following identities.
\begin{enumerate}
\item[(i)] $n_0>0\Leftrightarrow \mu\in I_n^0:$
\begin{eqnarray*}
&& H_\cG^*(\cA(\Si^\ell_i)^\pm_\mu;K)\cong H_\cG^*(\cN(\Si^\ell_i)^\pm_\mu;K)\cong
   H^*_{U(n)}(\ymU{\ell}{i}_\mu^\pm;K)\\
&\cong & H^*_{U(n_0)}(V_{ss}(\Si^\ell_i,P^{n_0,\pm(-1)^{n'i+k}} );K)\otimes
 \bigotimes_{j=1}^m H^*_{U(n_j)}(V_{ss}(\Si^{2\ell+i-1}_0, P^{n_j,k_j});K),
\end{eqnarray*}
\begin{eqnarray*}
&&  P^\cG_t(\cA(\Si^\ell_i)^\pm_\mu;K)= P^\cG_t(\cN(\Si^\ell_i)^\pm_\mu;K)=
   P^{U(n)}_t(\ymU{\ell}{i}_\mu^{\pm 1};K)\\
&=&P^{U(n_0)}_t(V_{ss}(\Si^\ell_i,P^{n_0,\pm(1)^{n'i+k}});K)\cdot
 \prod_{j=1}^m P^{U(n_j)}_t(V_{ss}(\Si^{2\ell+i-1}_0, P^{n_j,k_j});K).
\end{eqnarray*}
\item[(ii)]$n_0=0\Leftrightarrow \mu\in I_n^{i,\pm}:$
\begin{eqnarray*}
&& H_\cG^*(\cA(\Si^\ell_i)_\mu;K)\cong H_\cG^*(\cN(\Si^\ell_i)_\mu;K)\cong
H^*_{U(n)}(\ymU{\ell}{i}_\mu;K)\\
&\cong&\bigotimes_{j=1}^m H^*_{U(n_j)}(V_{ss}(\Si^{2\ell+i-1}_0, P^{n_j,k_j});K),\\
&&  P^\cG_t(\cA(\Si^\ell_i)_\mu;K)=  P^\cG_t(\cN(\Si^\ell_i)^+_\mu;K)=
    P^{U(n)}_t(\ymU{\ell}{i}_\mu;K)\\
&= &\prod_{j=1}^m P^{U(n_j)}_t(V_{ss}(\Si^{2\ell+i-1}_0, P^{n_j,k_j});K).
\end{eqnarray*}
\end{enumerate}

\end{thm}

\begin{ex} $n=2$.
For $\Si^\ell_1$,
$$
I_n^{1,+} =\{(0,0)\} \cup \{(2r-1,1-2r)\mid r\in \bZ_{>0}\},\quad
I_n^{1,-} =\{(0,0)\} \cup \{ (2r,-2r)\mid r\in \bZ_{>0}\}.
$$
For $\Si^\ell_2$,
$$
I_n^{2,+} =\{(0,0)\} \cup \{(2r,-2r)\mid r\in \bZ_{>0}\},\quad
I_n^{2,-} =\{(0,0)\} \cup \{ (2r-1,1-2r)\mid r\in \bZ_{>0}\}.
$$
In other words,  the $\cG$-equivariant strata of $\cA(P^{2,+})$(resp. $\cA(P^{2,-}$)
are $\{\cA(\Si^\ell_i)_\mu \mid \mu\in I^{i,+}_2\ \mbox{(resp. $I^{i,-}_2$)}\}$.
The codimension of each stratum is
$$
d_{2r,-2r}= 4r+ 2\ell+i-2,\quad
d_{2r-1,1-2r} = 4r+ 2\ell+i-4.
$$
The equivariant Poincar\'{e} polynomial for stratum $\mu=(r,-r)$ is
\begin{eqnarray*}
P^\cG_t\left(\cA(\Si^\ell_i)_{r,-r}\right)
&=& P^{U(2)}_t\left(X^{\ell,i}_{\mathrm{YM}}(U(2))_{r,-r}\right)
= P^{U(1)}_t\left(X^{2\ell+i-1,0}_{\mathrm{YM}}(U(1))_{r}\right) \\
&=& P^{U(1)}_t(U(1)^{2(2\ell+i-1)}) = \frac{(1+t)^{2(2\ell+i-1)}}{1-t^2}.
\end{eqnarray*}
\end{ex}

\begin{ex} $n=3$. Since $n$ is odd, $n_0=1$ or $3$. Thus,
$$
I_n^0=\{(0,0,0)\}\cup \{ (r,0,-r)\mid r\in\bZ_{>0} \}.
$$
The $\cG$-equivariant strata of $\cA(P^{3,\pm})$
are $\{\cA^\pm_\mu\mid \mu\in I^0_3\}$. The codimension of each stratum is
$$
d_{r,0,r}=4r+3(2\ell+i-2)
$$

The equivariant Poincar\'{e} polynomial for stratum $\mu=(r,0,-r)$ is
\begin{eqnarray*}
&& P^{U(3)}_t\left(X^{\ell,i}_{\mathrm{YM}}(U(3))_{r,0,-r}^{\pm 1}\right) =
P^{U(1)}_t\left(X^{\ell,i}_{\mathrm{flat}}(U(1))^{\pm(-1)^{i+r}}\right)
P^{U(1)}_t\left(X^{2\ell+i-1,0}_{\mathrm{YM}}(U(1))_{r}\right) \\
&&= P^{U(1)}_t(U(1)^{2\ell+i-1}) P^{U(1)}_t(U(1)^{2(2\ell+i-1)})
= \frac{(1+t)^{3(2\ell+i-1)}}{(1-t^2)^2}.
\end{eqnarray*}

\end{ex}

\end{document}